\documentclass{amsart}[10pt,draft]
\usepackage[colorlinks=true,          
            linkcolor=blue,
            citecolor=red,
            urlcolor=blue]{hyperref}
\usepackage{amssymb,amsfonts,amsmath,hyperref,mathrsfs,latexsym,stmaryrd,graphicx, manfnt}

\DeclareMathAlphabet{\mathpzc}{OT1}{pzc}{m}{it}

 \newtheorem{thm}{Theorem}[section]

 \newtheorem{Prop}{Proposition}[section]

    \theoremstyle{remark}\newtheorem{rmk}{Remark}[section]
    \theoremstyle{definition} \newtheorem{defn}{Definition}[section]
\input xy
\xyoption{all}

\newcommand {\Art}{\textrm{Art}_\CC}

\newcommand{\sets}{\textrm{Sets}}

\newcommand {\CC}{\mathbb{C}}

\newcommand {\GG}{\mathbb{G}}

\newcommand{\NN}{\mathbb{N}}
\newcommand {\PP}{\mathbb{P}}
\newcommand{\QQ}{\mathbb{Q}}
\newcommand {\RR}{\mathbb{R}}

\newcommand {\HH}{\mathbb{H}}

\newcommand{\ZZ}{\mathbb{Z}}

\newcommand {\ba}{{\bf a}}

\newcommand {\bg}{{\bf g}}

\newcommand {\bB}{{\bf B}}

\newcommand {\bG}{{\bf G}}

\newcommand {\bO}{{\bf O}}

\newcommand {\bR}{{\bf R}}

\newcommand {\bU}{{\bf U}}
\newcommand {\bV}{{\bf V}}

\newcommand{\cB}{\mathcal{B}}
\newcommand{\cC}{\mathcal{C}}

\newcommand{\cF}{\mathcal{F}}
\newcommand{\cG}{\mathcal{G}}
\newcommand{\cH}{\mathcal{H}}

\newcommand{\cL}{\mathcal{L}}

\newcommand {\cO}{\mathcal{O}}
\newcommand{\cP}{\mathcal{P}}
\newcommand{\cQ}{\mathcal{Q}}

\newcommand{\cX}{\mathcal{X}}

\newcommand{\cR}{\mathcal{R}}
\newcommand{\cS}{\mathcal{S}}
\newcommand{\cT}{\mathcal{T}}
\newcommand{\cU}{\mathcal{U}}

\newcommand{\cY}{\mathcal{Y}}

\newcommand{\scB}{\mathscr{B}}
\newcommand{\scC}{\mathscr{C}}

\newcommand{\scF}{\mathscr{F}}

\newcommand{\scK}{\mathscr{K}}
\newcommand{\scL}{\mathscr{L}}

\newcommand{\scP}{\mathscr{P}}
\newcommand{\scQ}{\mathscr{Q}}

\newcommand{\scZ}{\scZ}

\newcommand{\fb}{\mathfrak{b}}

\newcommand{\fg}{\mathfrak{g}}
\newcommand{\fh}{\mathfrak{h}}

\newcommand{\fl}{\mathfrak{l}}

\newcommand{\fs}{\mathfrak{s}}
\newcommand{\ft}{\mathfrak{t}}

\newcommand{\fA}{\mathfrak{A}}

\newcommand{\fD}{\mathfrak{D}}

\newcommand {\fU}{\mathfrak{U}}

\newcommand {\dbar}{\overline{\partial}}

\newcommand {\shom}{\textrm{\underline{Hom}}}
\newcommand{\mhom}{\textrm{Hom}}
\newcommand {\send}{\underline{End} }
\newcommand {\mend}{\textrm{End}}
\newcommand {\misom}{\textrm{Isom}}

\newcommand {\ad}{\textrm{ad} }

\newcommand{\spec}{\textrm{Spec }}
\newcommand{\sspec}{\underline{\textrm{Spec}}}
\newcommand {\cok}{\textrm{coker }}
\newcommand{\tot}{\textrm{tot }}

\newcommand{\ctimes}{\otimes_\CC}
\newcommand{\ztimes}{\otimes_\ZZ}

\newcommand  {\eps}{\varepsilon}

\newcommand {\io}{\iota}

\newcommand{\hookr}{\hookrightarrow}
\newcommand{\sym}{\textrm{Sym}}

\newcommand{\pr}{\textrm{pr}}
\newcommand{\rk}{\textrm{rk }}
\newcommand{\ram}{\textrm{Ram}}
\newcommand{\bra}{\textrm{Bra}}
\newcommand{\Higgs}{{\bf Higgs}}
 \newcommand{\Bun}{{\bf Bun}}

\newcommand{\Prym}{{\bf Prym}}

\newcommand{\crts}{\sf{coroot}}
\newcommand{\cwts}{\sf{coweight}}

\newcommand{\cchr}{{\sf cochar}}
\newcommand{\Aut}{\textrm{Aut}}

\newcommand {\ev}{\textrm{ev}}

\newcommand{\lestwo}[9]{
\xymatrix{     
 0 \ar[r] & {#1} \ar[r]  &  {#2} \ar[r]  &  {#3} 
\ar@{->}`r/10pt[d] `[l] `^dl[dlll]  `^dr/10pt[dll]    [dll] \\
 &  {#4} \ar[r] & {#5} \ar[r] & {#6} 
\ar@{->}`r/10pt[d] `[l] `^dl[dlll]  `^dr/10pt[dll]    [dll] \\
 & {#7} \ar[r]  & {#8} \ar[r] & {#9} }
}

\include{biblio}
\title{ Meromorphic Higgs bundles  And Related Geometries}
\author{Peter Dalakov}
\address{Institute of Mathematics and Informatics, Bulgarian Academy of Sciences, acad.G.Bontchev str., bl.8, 1113 Sofia Bulgaria}
\date{\today} 
\subjclass{14D20, 14C30, 14D15, 17B70, 53C07}
\newtheorem*{thma}{Theorem A}
\newtheorem{proposition}[thm]{Proposition}
\theoremstyle{definition}

\begin{document}

\begin{abstract}
The present note is mostly a survey on the generalised Hitchin integrable system and 
moduli spaces of meromorphic Higgs bundles. We also fill minor gaps in the existing literature,
outline a calculation of the infinitesimal period map and review briefly some related geometries.

\smallskip
\noindent \textbf{Keywords: Hitchin system, cameral covers, Donagi-Markman cubic, meromorphic Higgs bundles, integrable systems}
\end{abstract}

\maketitle 
\setcounter{tocdepth}{1}
\tableofcontents 

\section{Introduction}
    \subsection{Integrable systems and complex geometry}
Many moduli spaces arising  in complex-algebraic or analytic geometry   carry a symplectic or Poisson structure.
The  spaces  considered in this survey are no exception.  
Let $G$ be a simple complex Lie group,  $X$  a compact  Riemann surface with canonical bundle $K_X=\Omega^1_X$
and $D$ a sufficiently positive effective divisor on $X$.
Our  exposition is built
around the study of meromorphic, i.e.,  $K_X(D)$-valued, $G$-Higgs bundles on $X$ (Definition 2.1.) and their 
coarse moduli spaces 
$\Higgs_{G,D}$.
These  spaces come with  the additional structure of an algebraic completely integrable Hamiltonian system (ACIHS), known as the
generalised (or ramified) Hitchin system.

Completely integrable Hamiltonian systems  have long been an  object of interest for both mathematicians and physicists. 
The last thirty years have brought  significant advances  in the study of their  algebraic (and holomorphic) counterparts.
This was stimulated by the development of new methods in abelian and non-abelian Hodge theory, complex dynamics and holomorphic symplectic geometry, 
Yang--Mills and Seiberg--Witten theories, and of course, the quest for understanding mirror symmetry in its various incarnations.

The key difference between real and algebraic (or holomorphic) integrable systems is that abelian varieties (and complex tori) 
have moduli. Hence,  after the removal of singular fibres,
 the structure morphism of the ACIHS  is  a $C^\infty$  torus fibration, which usually fails to be holomorphically  locally trivial.
It is thus important to understand the corresponding period map, or, less ambitiously, the differential of the latter.

    \subsection{Contents of the paper}
We begin with  a  discussion of the moduli spaces 
$\Higgs_{G,D}$  (\S \ref{higgs_moduli}) 
and their Poisson geometry (\S \ref{symplectic}). Then in \S \ref{cameral}
we discuss  cameral covers and the Hitchin map.
 Most of the results in these introductory sections  are
standard and based on \cite{bottacin}, \cite{markman_thesis}, \cite{donagi_markman} and \cite{markman_sw}. 
There are, however, a number of well-known (and used)  extensions of results of Bottacin and Markman, for which we have not been able
to locate a proper reference. For these we have included partial proofs, wherever appropriate.

One of our goals in this note is to outline a calculation 
    of the infinitesimal period map of the
 generalised (ramified) Hitchin system.
This is done in \S  \ref{cubic}, and we refer to  \cite{ugo_peter_cubic} for more details. In short, our main result 
in \S \ref{cubic} is that the Balduzzi--Pantev formula (\cite{balduzzi}, \cite{ddp}) holds along the maximal rank symplectic leaves 
of the generalised  Hitchin system.

Admittedly, the ramified Hitchin system may seem very special, but we recall the
folklore statement that  all known 
ACIHS arise as  special case of Hitchin's.
 Some well-known examples
are geodesic flows on ellipsoids (Jacobi--Moser--Mumford system), KP elliptic solitons, Calogero--Moser  and elliptic Sklyanin systems. 
While some of these systems arise as complexifications of real CIHS, in general such complexifications
do not give rise to ACIHS, since
real Liouville tori need not ``complexify well''.
We direct the interested  reader to the wonderful surveys \cite{donagi_markman} and \cite{markman_sw} for
a detailed discussion and  examples.

Apart from Higgs bundles, cameral covers and Prym varieties, there
are several other geometric structures related to the space $\Higgs_{G,D}$:    special K\"ahler geometry, several flavours of
 Hodge theory, $tt^\ast$-geometry and Frobenius-like structures, to name a few.
We devote our final section \S \ref{survey}  to a  very brief literature review and  discussion of some of  these  structures.
    \subsection{Conventions and notation}
In \S \S \  \ref{higgs_moduli}, \ref{symplectic}, \ref{cameral} we alternate between the holomorphic and the algebraic viewpoint
and emphasise the differences, whenever deemed important.  For the proof of the main theorem in \S \ref{cubic}
we work in the holomorphic category. We  fix the following two types of ingredients:
    \begin{enumerate}
	  \item Geometric data:
	    \begin{itemize}
		    \item a smooth, compact, connected Riemann surface $X$ of genus $g\geq 0$
		    \item a  divisor $D\geq 0$ on $X$, with
  $K_X(D)^2$  very  ample 
	    \end{itemize}
	  \item Lie-theoretic data:
            \begin{itemize}
		    \item  a  simple complex Lie group $G$
		    \item   Cartan and Borel subgroups   $T\subset B\subset G$.
	    \end{itemize}
      \end{enumerate}
We denote by $Z$ or $Z(G)$ the centre of $G$.
The twist of the canonical bundle of $X$ by
$\cO_X(D)$ will be denoted by $L:= K_X(D)$. 
We shall  also  use the following -- mostly standard -- Lie-theoretic notation. 
The Lie algebras of the Cartan and Borel subgroups  will be
denoted, respectively, as $\ft\subset \fb\subset \fg$, while
$\cR^+\subset \cR\subset \ft^\vee$ will denote the (positive) roots.
 We let
$W=N_G(T)/T$ stand for   the abstract  Weyl group, which will be identified with its
embeddings in $GL(\ft^\vee)$ and $GL(\ft)$. Finally, let  $l=\rk \fg=\dim \ft$ be the rank of $G$, and 
$d_i$ ($1\leq i\leq l$)  the degrees of (any choice of) basic $G$-invariant polynomials on $\fg$.
For some of the calculations we will also use  a fixed choice of generators $\{I_i\}$ of $\CC[\fg]^G$.
We also fix an Ad-invariant symmetric bilinear form $\textrm{Tr}$ on $\fg$.

To these data one can associate two (closely related) families of abelian  torsors  parametrised by the \emph{Hitchin base}
 $\cB=H^0(X,\ft\ctimes L/W)\simeq H^0(X, \bigoplus_i   L^{d_i})$:
      \begin{itemize}
       \item a certain moduli space of meromorphic Higgs bundles on $X$
       \item a family of generalised Prym varieties for a family of (branched) $W$-Galois covers of $X$.
      \end{itemize}

Both (have connected components which) are  ACIHS in the Poisson sense.
The first family is  known as the ``generalised'' or ``ramified'' Hitchin system (with singular fibres removed).
The  second family is the ``abelianisation'' of the first one, and is (locally on the base) isomorphic it.  
While globally different, they  have the same infinitesimal period map, and  we shall use the second family
for our main Kodaira--Spencer calculation.

We remark that the Hitchin base $\cB$ depends  on  $G$, but  only via $\fg$, and we write
$\cB_\fg$ whenever it is important  to emphasise this dependence. 
There are certain  loci in $\cB$ for which we  introduce  special notation: the Zariski-open locus $\scB\subset \cB$
of generic cameral covers, $\cB_0\subset \cB$ for the locus (\ref{B0}) of pluri-differentials vanishing along $D$,
and $\bB\subset \scB$ for the base (\ref{boldB}) of the integrable system, obtained by restricting the Hitchin map to a maximal rank
symplectic leaf.

	\subsection{Acknowledgments} This note is a greatly expanded version of  my talk at the workshop on
 \emph{Instanton Counting: Moduli Spaces, Representation Theory and Integrable Systems}, which took place June 16-20, 2014 in Leiden, the Netherlands.
I would like to thank the organisers of the workshop: 
U.Bruzzo, D.Markushevich, V.Rubtsov, F.Sala and S.Shadrin, as well the staff of
 the Lorenz center for giving me the opportnuity
to speak and  for creating a  wonderful atmosphere of hospitality.

\section{Meromorphic $G$-Higgs Bundles}\label{higgs_moduli}
In this section we introduce our main objects of study: $G$-Higgs bundles on $X$ with values in a vector bundle. 
Next, we discuss the main global properties of the
 coarse moduli space of $K_X(D)$-valued $G$-Higgs bundles. 
Finally, we study in more detail the locus in the moduli space, corresponding to Higgs bundles whose
underlying principal bundle is regularly stable.

      \subsection{$L$-valued $G$-Higgs bundles}
	    \subsubsection{The definition}
Higgs bundles come in various flavours, and 
we begin by reviewing
one of the simplest variants of this notion:
 a Higgs bundle with values in a  vector  bundle.
    \begin{defn}
Let  $V$ be a holomorphic (algebraic) vector bundle on $X$.
 A holomorphic (algebraic)
 \emph{ $V$-valued $G$-Higgs bundle} on $X$ is  a pair $(P,\theta)$, consisting of
   a holomorphic (algebraic) principal $G$-bundle $P\to X$ and a section $\theta\in H^0(X, \ad P\otimes V)$,
called \emph{a Higgs field}.
    \end{defn}
A $V$-valued Higgs bundle is also called a Higgs bundle with $V$-coefficients, the rationale being
that for classical $G$, the Higgs field can be represented  in a local trivialisation by a matrix with coefficients
in $V$.
 This  terminology follows
\cite{don-gaits}, where the authors separate the r\^{o}le of the coefficient object from that
of the abstract (principal) Higgs bundle, see Definition 2.2, \emph{ibid.} 
For Higgs bundles with other coefficients, e.g., abelian fibrations, see Part V, \emph{ibid.}.
We remark that $V$-valued Higgs bundles,  especially when $V$ is a line bundle, are sometimes called $V$-twisted Higgs bundles. We
 avoid this terminology as potentially conflicting with the notion of a twisted Higgs bundle, understood as
a Higgs field on  a twisted bundle, the latter being  a
 bundle on a $\mu_n$-banded gerbe over $X$.

Our coefficient bundle $V$ will be,  almost exclusively, the \emph{line} bundle  $L= K_X(D)$. The only exception to this
is section \ref{survey}, where $V=T^\vee_\bB$, the cotangent bundle to a complex manifold, possibly non-compact and of
dimension  greater than one.
In the case $D=0$ our objects become 
\emph{$K_X$-valued $G$-Higgs bundles}, which are the objects  originally introduced by Hitchin  in \cite{hitchin_sd}
(for $G=SL_2(\CC)$, $PGL_2(\CC)$) and \cite{hitchin_sb} (for classical $G$).

An  isomorphism $(P,\theta)\simeq (Q,\phi)$ between  Higgs bundles  is  an isomorphism $f:P\simeq Q$
 between the $G$-bundles, which preserves the Higgs fields, i.e., $\ad f(\theta)=\phi$.

The above notion of Higgs bundles is clearly functorial in the structure group.  If $(P,\theta)$ 
is  an $L$-valued $G$-Higgs bundle, and
  $\rho:G\to H$  a group
homomorphism, then extension of  structure group gives an $H$-bundle $P'= P\times ^G H$.
Moreover, we have a homomorphism $\ad\circ \rho: G\to \Aut(\fh)$, and the associated bundle
$P(\fh)= P\times^{\ad\circ \rho}\fh$ is isomorphic to $\ad P'$. Thus $\theta$ gives rise,
via the homomorphism $\fg\to \fh$, to a Higgs field  $\theta'$ on $P'$.
	    \subsubsection{Families}
If $S$ is a complex manifold or complex space (respectively, an algebraic variety or scheme over $\CC$), a \emph{family} of 
$L$-valued $G$-Higgs bundles on $X$, parametrised by $S$ is a pair $(\scP, \Theta)$, where $\scP\to S\times X$ is a holomorphic (respectively,  algebraic) principal
$G$-bundle and  $\Theta$ is a section (over $S\times X$) of $\ad\scP\otimes p_X^\ast L$.
 Equivalently, we think of $\Theta$ as being a section (over $S$)
 of $p_{S \ast} (\ad\scP\otimes p_X^\ast L)$, where $p_S=\textrm{pr}_1$ stands for the canonical projection.
 In the case of a pointed   base $(S,o)$,  $o\in S$, we call a family $(\scP,\Theta)$ a 
\emph{deformation} of $(P,\theta)$ if it is equipped with an isomorphism $(\scP_o,\Theta_o)\simeq(P,\theta)$.
There is also an obvious notion of isomorphism of families (and deformations). 

In the sequel we suppress the distinction between the algebraic and the analytic case unless there is a danger of confusion.
For vector bundles on curves this can be justified by the GAGA principle. In higher dimensions (e.g., for families)
one should keep in  mind that algebraic $G$-bundles are assumed to be isotrivial, i.e., trivial in \'etale topology,
rather than in the Zariski topology, see \cite{serre_fibres_algebriques}. In case when $G=GL_n(\CC)$ (which we exclude)
or when $G$ is reductive and the base is a curve, one \emph{can} indeed use the Zariski topology, by a result
of Springer (\cite{steinberg_regular}). Moreover, by a theorem of Drinfeld and Simpson (\cite{drinfeld_simpson})
$G$-bundles on $S\times X$ are locally trivial in the product of the \'etale topology on $S$ and Zariski topology on
$X$.

The following elementary example of a family of Higgs bundles will be needed in what follows. 
Consider $S=H^0(X, \ad P\otimes L)$, where   $P\to X$ is a  $G$-bundle. 
Then the trivial family of $G$-bundles $\scP= p_X^\ast P\to S\times X$ can be augmented to a
 family of Higgs bundles.
Indeed,  
\begin{equation}\label{tautological}
p_{S \ast}\left(\ad\scP\otimes p_X^\ast L\right) = \cO_S\ctimes H^0(X,\ad P\otimes L)= T_S,
 \end{equation}
 and
we take
 $\Theta\in H^0(S, p_{S \ast} (\ad\scP\otimes p_X^\ast L)) = H^0(S, T_{S})$  
  to be the tautological section of $T_S$, i.e., the Euler vector field on $S$.

\vspace{5pt}

	  \subsection{Moduli Spaces}
	      \subsubsection{Principal Bundles}
Given a $G$-bundle $\pi: P\to X$ and a   closed algebraic subgroup  $R\subset G$, we obtain an
associated $G/R$-bundle $\pi_R: P\times^G (G/R)= P/R\to X$ and a
 principal $R$-bundle $P\to P/R$, which we also denote by $P_R$, to avoid confusion.
This relies on Proposition 3 of 
\cite{serre_fibres_algebriques}, stating  that $G\to G/R$ is a principal
$R$-bundle (in the \'etale topology). In the  analytic category,
for $R\subset G$  a closed complex Lie subgroup,
this  follows from a theorem of Chevalley about existence of local analytic  sections of the canonical projection.
We then identify the set of $R$-reductions of
$P$ with the set of sections $\Gamma(X,P/R)$  in the usual way: a section $\sigma: X\to P/R$
 gives rise to  an $R$-reduction $\sigma^\ast P_R= X\times_{P/R}P\subset P$.

 Following
Ramanathan (\cite{ramanathan}, Definition 1.1), we say that
a principal bundle $P\to X$ is \emph{stable} (respectively, \emph{semi-stable}) if, for every maximal parabolic subgroup
$H\subset G$, and every  $H$-reduction $\sigma: X\to P/H$, $\deg \sigma^\ast T_{\pi_H} >0$ 
(resp. $\deg \sigma^\ast T_{\pi_H} \geq 0$). Here
\[
T_{\pi_H}= \ker (d\pi_H)  \subset T_{P/H}
\]
stands for the relative tangent 
bundle  of the morphism $\pi_H$ and is nothing but $P_H(\fg/\fh)= \left(P\times \fg/\fh\right)/H$. Equivalently (Lemma 2.1, \emph{ibid.}),  $P$ is (semi-)stable if
for any  reduction $\sigma:X\to P/H$ to a parabolic subgroup $H\subset G$, and any dominant character
$\chi$ on $H$, one has $\deg \sigma^\ast\left( E_H\times^\chi \CC^\times \right)<0$ (respectively, $\leq 0$).

Ramanathan constructed a (coarse) moduli space $\Bun_G$ of S-equivalence classes of semi-stable $G$-bundles
(or isomorphism classes of poly-stable $G$-bundles).
By Theorem 4.3 (\cite{ramanathan}), $\Bun_G$ is a normal Hausdorff analytic space, with connected components
$\Bun_{G,c}$, indexed by the topological type $c\in \pi_1(G)$ of the bundle:
\[
 \Bun_G= \coprod_{c\in\pi_1(G)} \Bun_{G,c}.
\]
Ramanathan also constructed the moduli space in the algebraic category, identifying $\Bun_{G,c}$  as the analytification of a normal projective algebraic variety
  (\cite{rama1}, \cite{rama}).
For a Tannakian construction of the moduli space, see \cite{balaji_seshadri_1}.

By \cite{ramanathan}, Proposition 3.2,  if $P$ is stable, then $H^0(X,\ad P)=0$ and $\Aut P$ is finite.
 More generally, in  the case of a reductive -- rather than  simple -- group $G$ one has $H^0(X,\ad P)=\textrm{Lie }Z(G)$, \emph{ibid.}.
Since, if $g\geq 2$, every topological $G$-bundle admits some structure of a stable holomorphic $G$-bundle (Remark  5.3, \emph{ibid.}),
we see by a Hirzebruch--Riemann--Roch calculation  that
\[
\dim \Bun_G= \dim G\cdot  (g-1).
\]
For a discussion of low-genus cases, see
\cite{Grothendieck_P1}, \cite{atiyah_elliptic}, \cite{tu_elliptic}, \cite{laszlo_ell}, \cite{friedman_morgan_witten}.

 The infinitesimal deformations of a semi-stable $G$-bundle $P$ are parametrised by
$H^1(X,\ad P)$, and, by Luna's \'etale slice theorem, the GIT quotient
\[
 H^1(X,\ad P)\sslash \Aut (P)
\]
is isomorphic to an \'etale neighbourhood of  $[P]\in \Bun_{G}$.  
There is a natural inclusion  $Z(G)\subset \Aut(P)$, and, for
  $g\geq 2$, the smooth locus of the open subvariety of stable bundles $\Bun^{st}_G\subset \Bun_G$ consists of the regularly stable bundles, i.e.,
those which satisfy $\Aut P= Z(G)$, see \cite{ramanathan} or \cite{biswas_hoffmann_torelli_bun}, Proposition 2.3.
		\subsubsection{Higgs Bundles}
Ramanathan's definition of (semi-)stability also makes sense for
  Higgs bundles, provided one considers only parabolic reductions which ``preserve the Higgs field'', in the following sense.
Given a closed subgroup $R\subset G$ and an  $R$-reduction $\sigma: X\to P/R$, we have a natural  projection
\[
 \xymatrix@1{ \Pi_\sigma:  &\ad P\ar[r]&\ \frac{\ad P}{\ad\ \sigma^\ast P_R} }.
\]
If $\theta$ is a Higgs field on $P$, we say that an $R$-reduction $\sigma$ of $P$
 is a \emph{Higgs reduction} of $(P,\theta)$
if $\theta\in\ker (\Pi_\sigma\otimes id)$.
If $\sigma$ is a Higgs reduction, then the Higgs field 
 $\theta$ on the $G$-bundle $P$ induces  a Higgs field on the $R$-bundle $\sigma^\ast P_R$.

In this way, the choice of Higgs field singles out a class of Higgs reductions among all $R$-reductions.
 This class can be conveniently described in   the approach from \cite{ugo_bea_ss}, Definition 3.5. The projection
$\fg\to \fg/\mathfrak{r}$ induces a bundle homomorphism 
\[
\eta: \pi_R^\ast \ad P=P_R(\fg)\longrightarrow T_{\pi_R}= P_R(\fg/\mathfrak{r}),
\]
and thus $\theta$ gives rise to a section
\[
 \eta(\pi_R^\ast \theta)\in H^0(P/R, T_{\pi_R}\otimes \pi_R^\ast K_X(-D))\subset H^0(P/R, T_{\pi_R}\otimes \Omega^1_{P/R}(-D)).
\]
The vanishing locus of this section determines a closed  subscheme in $P/R$,  \emph{the scheme of Higgs reductions} of $(P,\theta)$. A reduction
$\sigma: X\to P/R$ is a Higgs reduction precisely when its image  is contained in the scheme of Higgs reductions. This scheme
turns out to play an important r\^ole  in studying $\Omega^1_X$-valued $G$-Higgs bundles on smooth projective varieties, but can be, in general, rather 
singular.

 We say that $(P,\theta)$ 
is \emph{(semi-)stable} if, for any  Higgs reduction $\sigma:X\to P/H$ to a maximal parabolic $H\subset G$,
$\deg T_{\pi_H}>0$ (respectively, $\deg T_{\pi_H}\geq 0$).
Suitably modifying Ramanathan's construction, 
one obtains a quasi-projective
coarse  moduli space $\Higgs_{G,D}$ of S-equivalence classes of semi-stable Higgs bundles.
When $D=0$, $\Higgs_{G,0}$ is  in fact 
 a partial compactification of $T^\vee\Bun^{sm}_G$ (\cite{hitchin_sd}, \cite{hitchin_sb}).
Moreover, when $D=0$, it is known from \cite[Lemma 4.2]{don-pan} (see  \cite{oscar_andre_compts} for a different proof) 
that the connected components of the moduli space are indexed by $\pi_1(G)$, i.e., by the topological type of the 
 $G$-bundle, underlying  a Higgs bundle. 
This is expected to hold for
arbitrary $D>0$ (whenever  the moduli space is non-empty), but there does not seem to exist a published statement to this effect.
However, for each $c\in \pi_1(G)$ there exists an irreducible connected component, $\Higgs_{G,D,c}$, characterised by the fact that it
 contains generic cameral covers, see \ref{cameral_covers} and \cite{donagi_markman}, definition 4.9.
 As we shall see \S \ref{cameral}, 
there is a morphism (the Hitchin map) from
$\Higgs_{G,D}$ to a vector space $\cB$, and a strict subvariety, $\Delta\subset \cB$, such that the connected components of the fibres of
$\left. \Higgs_{G,D}\right|_{\cB-\Delta}\to \cB-\Delta$ are isomorphic to abelian varieties. These connected components of the Hitchin fibre
are contained in the respective connected components $\Higgs_{G,D,c}$. It seems presently unknown whether there exist connected components of
$\Higgs_{G,D}$, lying entirely over the discriminant locus $\Delta\subset \cB$, but  the arguments in \cite{don-pan}, Lemma 4.2 seem to indicate that
this is not the case.

For a discussion of S-equivalence, Harder--Narasimhan and Jordan--H\"older reductions in the case $D=0$ (but possibly over
higher-dimensional base), see \cite{arijit_parthasarathi_hn}, \cite{bea_jh}, \cite{ugo_bea_ss}.

 \'{E}tale locally near
$[(P,\theta)]$  the moduli space $\Higgs_{G,D,c}$ is isomorphic to
\[
  \HH^1 (\scK^\bullet_{(P,\theta)})\sslash \Aut (P,\theta),
\]
where $\scK^\bullet_{(P,\theta)}$ is the deformation complex (\ref{br_mero}).
A stable pair $(P,\theta)$ represents a smooth (regular) point  in $\Higgs^{st}_{G,D=0}$     precisely when it  is
regularly stable, i.e.,   $\textrm{Aut}(P,\theta)=Z(G)$, see 
\cite{Biswas-Ramanan} Theorem 3.1 or
\cite{faltings}. Then by a Hirzebruch--Riemann--Roch calculation
(see section \ref{def_th}, equation (\ref{dim}) ) one obtains
\begin{equation}\label{higgs_dim}
\dim \Higgs_{G,D}= \dim G\deg K_X(D).
\end{equation}
 In particular, for $D=0$ (which, with our assumptions, implies
$g\geq 2$) we have $\dim \Higgs_{G,D=0}= 2\dim \Bun_G$. The moduli space $\Higgs_{G,D=0}$ is  normal,  with orbifold singularities at worst.

One can also construct the moduli space in the algebraic category,
 following a  version of either Simpson's (\cite{moduli2} \S 4) or Nitsure's 
( \cite{nitsure} \S 5)
construction. The former deals with $D=0$, while the latter with $G=GL_n$.
We should mention here  that C.Simpson's notion of semi-stability is not always equivalent to Ramanathan's, see
\cite{ugo_bea_ss}, Remark 4.6. 
For a purely algebraic, GIT-free construction of the moduli space
in the case $D=0$, see \cite{faltings}.
If willing to work only with \emph{everywhere regular} Higgs fields, one can construct a moduli space via   the spectral correspondence 
(\cite{donagi_spectral_covers}, \S 5.4). 
We recall that  $\fg\backslash \fg^{reg}\subset \fg$ is of codimension three and, since we are considering only 
(line bundle-valued)
 Higgs bundles
on \emph{curves}, being everywhere regular is a reasonable restriction. 

As with principal bundles, one may prefer to fix a linear representation of $G$, and
work with vector bundles with extra structure. An $L$-valued
Higgs \emph{vector} bundle is a  pair $(E,\theta)$, consisting of 
 a vector bundle $E$ and a section $\theta\in H^0(X,\send E\otimes L)$.
 For such pairs one defines stability using
slope: $(E,\theta)$ is \emph{(semi-)stable}, if, for every subbundle $F\subset E$, 
satisfying
$\theta (F)\subset F\otimes L$, the inequality  $\mu (F)< \mu (E)$, respectively 
$\mu (F)\leq \mu (E)$, holds; see \cite{hitchin_sd}, \cite{nitsure}. When $\theta=0$, this reduces to Mumford's original notion of
stability of vector bundles.
We do not delve into a detailed comparison of the different notions of stability and the different constructions of moduli spaces, 
 mainly because we are going to work exclusively with
 generic (regularly stable) Higgs bundles. We do, however, discuss briefly the behaviour of stability under 
group homomorphisms and  compare  the stability of a Higgs bundle with the stability of its adjoint bundle.  In the next theorem,
we relax slightly our usual assumptions and allow \emph{reductive} structure groups.

\begin{thm}\label{stab_functoriality}
Let $X$ be a compact Riemann surface and $D$ an effective divisor on $X$, such that $H^0(X, K_X(D) )\neq (0)$. 
      \begin{enumerate}
       \item Let $V$ be a finite-dimensional $\CC$-vector space,  $G=GL(V)$,  and $P$ a (holomorphic) principal $G$-bundle, 
so that $\ad P=\send (P\times ^G V)= P\times ^G \mend (V)$.
Then   $(P,\theta)$ is a stable (semi-stable) Higgs bundle if and only if $(P\times^G V,\theta)$ is a stable (semi-stable) Higgs vector bundle.
       \item   Let $\phi:G\to H$ be a \emph{surjective} homomorphism
between reductive (complex) algebraic groups, such that $\ker \phi\subset Z(G)$. 
 Let $(P,\theta)$
be an $L$-valued $G$-Higgs bundle and $(P',\theta')$ the $L$-valued $H$-Higgs bundle, induced by $\phi$. 
Then $(P',\theta')$ is
stable (semi-stable) if and only if  $(P,\theta)$ is so.
      \item An $L$-valued $G$-Higgs bundle $(P,\theta)$ is semi-stable if and only if the adjoint Higgs (vector) bundle
$(\ad P,\ad\theta)$ is semi-stable. If $(\ad P,\ad \theta)$ is \emph{stable}, then so is
$(P,\theta)$. If $(P,\theta)$ is stable, $(\ad P,\ad\theta)$ need not be stable, but is polystable.

	\item A Higgs bundle $(P,\theta)$ is semi-stable if and only if for any representation
$\phi: G\to \Aut(V)$, such that
$\phi\left( Z_0(G)\right)\subset Z(\Aut(V))$ the associated Higgs vector bundle is semi-stable.

      \end{enumerate}

\end{thm}
\emph{Proof:}
Maximal parabolic subgroups $H\subset GL(V)$ 
consist of automorphisms of $V$, preserving a flag
$(0)\neq U\subsetneq V$, $V\simeq U\oplus V/U$, and $\sigma: X\to P/H$ is a Higgs reduction precisely when $\theta$ preserves the vector bundle 
$P(U)=P\times^G U\subset P(V)$.
Moreover, $\sigma^\ast T_{\pi_H}= \mhom (P(U), P(V/U))$, and since
\[
 \deg  \mhom (P(U), P(V/U)) = \rk P(U)\rk P(V/U)\left(\mu(P(V/U))-\mu(U)\right),
\]
statement (1) is proved.
 See also  \cite{arijit_parthasarathi_hn} Lemma 7,
 \cite{ramanathan} Lemma 3.3 and \cite{hyeon_murphy},  Corollary 1.

Statement (2) is proved as  Proposition 7.1 in \cite{ramanathan}.
The key point is that there is a one-to-one
correspondence 
between the parabolic reductions of $P$ and those of $P'$.
Indeed,  one sees that the diagram
\[
 \xymatrix{ 1\ar[r]& \ker\phi \ar[r]  & G\ar[r]^-{\phi} & H\ar[r]&1 \\
	    1\ar[r]&\ker\phi\ar@{=}[u]\ar[r]& R\ar[u]\ar[r]   & R'\ar[u]\ar[r]& 1  }
\]
induces, if $\ker\phi\subset Z(G)$, maps between the corresponding cohomology groups
\[
	\xymatrix{   
		  H^1(X,\ker\phi(\cO_X))\ar[r]           & H^1(X, G(\cO_X))\ar[r]       & H^1(X, H(\cO_X))\ar[r]        & H^2(X,\ker\phi (\cO_X))\\
		  H^1(X,\ker\phi(\cO_X))\ar[r]\ar@{=}[u] & H^1(X, R(\cO_X))\ar[r]\ar[u] & H^1(X, R'(\cO_X))\ar[r]\ar[u] & H^2(X,\ker\phi (\cO_X))\ar@{=}[u]},
\]
see  \cite{Grothendieck_kansas}, 5.7.11. Here $G(\cO_X)$ denotes the sheaf of germs of holomorphic maps from $X$ to $G$.
Since this correspondence preserves Higgs reductions, (2) follows.

The first part of
  statement (3) is proved as in the case $D=0$, for which we refer to \cite{biswas_anchouche_eh}, Lemma 4.7, 
 \cite{arijit_parthasarathi_hn} Proposition 12 and
\cite{ugo_bea_ss}, Lemma 4.3 (i). Notice that in these arguments one can
use statement (1) to pass from $G$ to $G^{ad}= G/Z(G)$.

The second part of statement (3) is proved by modifying the corresponding 
argument for principal bundles (e.g., Proposition 2 in \cite{hyeon_murphy}). Indeed, let $H\subset G$ be a maximal parabolic
and $\fh=\textrm{Lie }H$. The short exact sequence of vector spaces
\[
 \xymatrix@1 {0\ar[r] & \fh\ar[r]& \fg\ar[r] & \fg /\fh\ar[r]& 0}
\]
is a sequence of $H$-modules via the adjoint representation $H\hookr G\to \Aut (\fg)$. Twisting with it the
$H$-bundle $P_H = P\to P/H$ gives a sequence of vector bundles on $P/H$, which, when pulled back by an $H$-reduction
$\sigma: X\to P/H$, gives
\[
 \xymatrix@1{0\ar[r] & \sigma^\ast (P_H\times^ H \fh)\ar[r] & \ad P\ar[r] & \sigma^\ast T_{\pi_H}\ar[r] & 0  }.
\]
Suppose that $(\ad P,\ad \theta)$ is a stable Higgs (vector) bundle. If the above reduction $\sigma$ is a Higgs reduction, 
then $ \sigma^\ast (P_H\times^H \fh)\subset \ad P$ is preserved by $\ad\theta$ and hence, by stability, 
$\deg  \sigma^\ast (P_H\times ^H \fh) <\deg \ad P$. But $\deg \ad P=0$, since $G$ is reductive and a choice of $G$-invariant bilinear form on $\fg$ gives
an isomorphism $\ad P\simeq \ad P^\vee$. Hence $\deg \sigma^\ast T_{\pi_H}>0$.

For the final part of (3), see \cite{biswas_anchouche_eh}, Theorem 4.8. Examples of stable Higgs bundles
which are not ad-stable (but strictly ad-semistable) exist already for principal bundles, i.e., when $\theta=0$.
Moreover, such bundles  always exist if $\dim Z(G)>0$.

For part (4), see \cite{biswas_anchouche_azad_prb}, Lemma 1.3 and
 \cite{ugo_bea_ss}, Lemma 4.3 (ii).
\qed

      \subsection{Over the locus of regularly stable bundles}\label{reg_st}
By Theorem II.6 in \cite{faltings}, if $X$ is of genus $g\geq 2$, one has
$\Bun_G^{rs}\neq \varnothing$. Moreover,  if $g\geq 3$ or if $g\geq 2$ but $G\neq PGL_2$, then
the codimension of the complement of $\Bun_G^{rs}$ in $\Bun^{st}_G$ is at least two. 
In this section we   assume this to be the case, and
 consider  the Zariski open $\Higgs_{G,D}^o\subset \Higgs_{G,D}$, consisting of classes of pairs 
$(P,\theta)$, for which $[P]\in \Bun_G^{rs}$.  We sketch a direct construction of this locus 
as a vector bundle over $\Bun_G^{rs}$  and discuss the existence of Poincar\'e family.
Most of the  constructions in this section are
a natural   generalisation of \cite{bottacin}, \S\S 1, 3.

The short exact sequence
\[
 \xymatrix@1{ 1\ar[r] & Z(G)\ar[r] & G\ar[r] & G^{ad}\ar[r] & 1 }
\]
gives rise to an exact sequence of pointed sets
\[
\xymatrix{   	 H^1(X, Z(G))\ar[r]           & H^1(X, G(\cO_X))\ar[r]       & H^1(X, G^{ad}(\cO_X))\ar[r]        & H^2(X,Z(G))}
\]
and to a morphism $\pi: \Bun_G\to \Bun_{G^{ad}}$. On closed points the latter is given by
$\pi([P])= [P']$, where $P'= P\times^G G^{ad}= P/Z(G)$.
Let us fix a topological type $c\in \pi_1(G)$ and consider the restriction
\[
 \pi_c: \Bun^{rs}_{G,c}\longrightarrow \Bun^{rs}_{G^{ad},c'},
\]
where $c'\in \pi_1(G^{ad})$ is the image of $c$ under the injection
$\pi_1(G)\hookr \pi_1(G^{ad})$ induced by the covering space $G\twoheadrightarrow G^{ad}$.
 Notice that $\pi$ respects both stability (by
Theorem \ref{stab_functoriality}) and minimality of automorphisms.
By \cite{biswas_hoffmann_poincare}, Corollary 6.9
(see also  \cite{bal_bis_nag_news}, Theorem 1.1 if $c'=0$)
 there exists a universal  $G^{ad}$-bundle  $\scP'\to \Bun_{G^{ad},c'}^{rs}\times X$.
Its pullback $\pi_c'^{\ast} \scP'$, $\pi_c'=(\pi_c,1)$,  could be called an \emph{adjoint Poincar\'e bundle}, since
$\left. \scP'\right|_{\{E\}\times X} \simeq E/Z(G)$, for $[E]\in \Bun_{G,c}^{rs}$.

To incorporate Higgs fields, 
 consider the adjoint bundle of $\pi_c'^\ast\scP'$, i.e., the
vector bundle $\pi_c^\ast \ad \scP'= \ad \left(\pi_c^\ast\scP'\right)$ on $\Bun_{G,c}^{rs}\times X$.
By semi-continuity and Grauert's theorem, 
the quasi-coherent sheaf 
$\cF = p_{1 \ast}\left(\ad\left(\pi_c^\ast\scP'\right)\otimes p_2^\ast L\right)$
 is locally free of
finite rank, and its total space
$\tot\cF= \sspec\  \sym^\bullet \cF^\vee$ is a vector bundle on $\Bun^{rs}_{G,c}$. 
Considering, for any $[E]\in\Bun_{G,c}$,  the diagram
\[
 \xymatrix{ \{E\}\times X\ar@<-0.2ex>@{^{(}->}[r]_-{j_E'}\ar[d]_-{p_1} & \Bun^{rs}_{G,c}\times X\ar[d]^-{p_1}\\
	    \{E\}\ar@<-0.2ex>@{^{(}->}[r]_-{j_E} & \Bun^{rs}_{G,c}\\   }
\]
and using that $\ad E=\ad \left(E/Z(G)\right)$, we obtain a canonical identification between the fibre of
$\cF$ over $[E]$ and the vector space of $L$-valued Higgs fields on $E/Z(G)$:
\[
 \cF_{[E]}=j_E^\ast \cF = j_E^\ast p_{1 \ast}\left(\ad \pi_c^\ast\scP'\otimes p_2^\ast L\right) = p_{1\ast}\left(\ad E\otimes p_2^\ast L\right) = H^0(X, \ad E\otimes L).
\]
We thus have
identified  $\tot \cF$ with  $\Higgs_{G,D,c}^o$.
Notice that if $D=0$, we have $\cF\simeq \Omega^1_{\Bun^{rs}_{G,c}}$.

Finally, we turn to the question of
existence of Poincar\'e family of Higgs bundles on $\Higgs_{G,D,c}^o$.  By general arguments (Luna's \'etale slice theorem),
such a family always exists locally in the \'etale (or analytic) topology. Using some recent results of Biswas--Hoffmann
and Donagi--Pantev,  
we can say a bit more about the global or Zariski-local situation as well.

Given a non-empty  Z-open $U\subset \Higgs_{G,D,c}^o$,  an \emph{adjoint Poincar\'e family} of $L$-valued $G$-Higgs bundles on $U$
is a pair $(\scQ,\Theta)$, where $\scQ\to U\times X$ is an adjoint Poincar\'e family of $G$-bundles, i.e., a $G^{ad}$-bundle, satisfying 
$\left. \scQ\right|_{\{[E,\theta]\}\times X}\simeq E/Z(G)$, while  $\Theta\in H^0(U, p_{U\ast}\left(\ad \scQ\otimes p_X^\ast L\right))$
is a family of $L$-valued $G$-Higgs bundles over $U$, 
such that $\left. \Theta\right|_{\{[E,\theta]\}}=\theta$.

We now introduce an extra piece of notation, following  \cite{biswas_hoffmann_poincare}.
Consider the coroot, cocharacter and coweight lattices in $\ft$:
\[
 \crts_\fg\subset \cchr_G\subset \cwts_\fg.
\]
These lattices can be identified with $\mhom (\GG_m, T^{sc})$, $\mhom (\GG_m, T)$ and
$\mhom (\GG_m, T^{ad})$, respectively, where
\[
  T^{sc} \twoheadrightarrow T \twoheadrightarrow T^{ad}
\]
are maximal tori in $G^{sc}$, $G$ and $G^{ad}$, respectively.
Correspondingly, the fundamental groups of $G$ and $G^{ad}$
are
\[
 \pi_1 (G)=\frac{\cchr_G}{\crts_\fg}\subset \pi_1(G^{ad})=\frac{\cwts_\fg}{\crts_\fg}.
\]

As one can  see (\cite{biswas_hoffmann_poincare}, Lemma 6.2), any
even, $W$-invariant, integral, symmetric bilinear form on $\crts_\fg$ extends to a symmetric, $\QQ/\ZZ$-valued bilinear form
on $\pi_1(G^{ad})$. These extensions generate a cyclic group 
\[
\Psi(G^{ad})\subset \mhom (\pi_1(G^{ad})^{\otimes 2},\QQ/\ZZ),
\]
cf. Table 1 in \cite{biswas_hoffmann_poincare}. Consider the subgroup
\[
 \Psi'(G) = \left\{\left. b\in \Psi(G^{ad})\right|  \left. b\right|_{\pi_1(G)\times \pi_1(G)}=0 \right\}\subset \Psi(G^{ad})
\]
of bilinear forms, vanishing on $\pi_1(G)$.
With every element $c\in\pi_1(G)$ we associate an evaluation map
\[
 \ev_G^c: \Psi'(G)\to \mhom \left(\frac{\pi_1(G^{ad})}{\pi_1(G)}, \QQ/\ZZ\right),\ b\mapsto b(c,\ ).
\]
More generally,
Biswas and Hoffmann  define an analogue of $\Psi'(G)$ for an arbitrary reductive group $G$ (Definition 6.4, \emph{ibid.}). They  tie
the obstruction of the existence of Poincar\'e family on $\Bun_{G,c}^{rs}$ with the cokernel of $\ev_G^c$ (Theorem 6.8, \emph{ibid}.).
More concretely, the moduli stack of regularly stable bundles (of type $c$) is a $Z(G)$-banded gerbe over $\Bun_{G,c}^{rs}$, and the order
of its class can be expressed via the order of finite group $\cok\ev_G^c$. For a related result in the Higgs setting, 
see  Lemma 4.2 in \cite{don-pan}.

      \begin{thm}
       Let $G$ be a simple complex algebraic group, $c\in\pi_1(G)$ and $X$ a compact Riemann surface of genus  $g\geq 3$ (or $g\geq 2$ and
$G\neq PGL_2$). Then there exists an \emph{adjoint} $L$-valued Poincar\'e $G$-Higgs bundle $(\scQ,\Theta)$ over $\Higgs_{G,D,c}^o\times X= \tot \cF\times X$.
There exists a non-empty open subscheme $U\subset \Bun_{G,c}^{rs}$ and an $L$-valued Poincar\'e $G$-Higgs bundle $(\scP,\Theta)$ over
$\tot \cF_{U}\subset \Higgs_{G,D,c}^o$ if and only if $\cok \ev_G^{c}=0$. If such an open $U$ exists,
then the Poincar\'e family extends to all of $\Higgs_{G,D,c}^o$ if and only if
 $G=G^{ad}$.
      \end{thm}

\emph{Proof:}

 Consider the diagram
\begin{equation}\label{universal}
 \xymatrix{	\tot\cF\times X\ar[r]_-{p_\cF'}\ar[d]^-{p_1}\	 & \Bun_{G,c}^{rs}\times X \ar[r]_-{\pi_c'}\ar[d]^-{p_1}\ & \Bun_{G^{ad},c'}^{rs}\times X\ar[d]^-{p_1}\\
		\tot\cF\ar[r]_-{p_\cF}	                     \	 & \Bun_{G,c}^{rs}\ar[r]_-{\pi_c}			\  &\Bun_{G^{ad},c'}^{rs}\\
},
\end{equation}
where $p_{\cF}: \tot\cF\to \Bun_{G,c}^{rs}$ is the bundle projection, $p_\cF([P,\theta])=[P]$. 
We have already  seen that $\pi_c'^\ast\scP'$ is an adjoint Poincar\'e bundle. Now we pull it further back 
 and set $\scQ= (\pi_c'\circ p_\cF')^\ast \scP'\to \tot\cF\times X$.
To construct the Higgs field on it, recall that, 
as with any vector bundle,
 $p_\cF^\ast \cF$ carries a tautological section, $\lambda\in H^0(\tot \cF,p_\cF^\ast \cF)$.
But then
\[
 p_\cF^\ast p_{1\ast} \left( \pi_c'^\ast \ad\scP'\otimes p_2^\ast L\right)\simeq p_{1 \ast}\left((\pi_c\circ p_\cF)^\ast\ad\scP'\otimes p_2^\ast L\right)
\]
and we get the family of $G^{ad}$-Higgs bundles $(\scQ,\Theta)= \left((\pi_c'\circ p_\cF')^\ast \scP',\lambda\right)$, which is an adjoint Poincar\'e family.

Since $\ad \scQ=\ad \left(\scQ/Z(G)\right)$, 
whenever there exists a Poincar\'e $G$-bundle on $\Bun_{G,c}^{rs}\times X$ (or an open subscheme thereof), we can pull the latter
by $p_\cF'$ to a $G$-bundle $\scP\to \tot \cF\times X$ and obtain a Poincar\'e family $(\scP,\lambda)$. Then the conditions for the existence of a
regularly stable Poincar\'e $G$-bundle are given in  Corollary 6.7, 6.9 of \cite{biswas_hoffmann_poincare} and Remark 6.10, \emph{ibid.}.

\qed

\section{Poisson Geometry}\label{symplectic}
In this section we discuss the symplectic and Poisson aspects of the geometry of Higgs moduli.
As a means of motivation, we start with $K_X$-valued Higgs bundles and recall the construction of the symplectic
form on $\Higgs_{G,0,c}$. Next we discuss the deformation theory of Higgs bundles, following Biswas and Ramanan, 
and describe the corresponding deformation complex. A Poisson bivector on a variety  determines a morphism from
the cotangent to the tangent sheaf of the latter. In our context, tangent spaces are expressed as hypercohomology groups of
 complexes of sheaves. We review duality for hypercohomology in \ref{duality_sect} and define the Poisson bivector
in \ref{poisson_sect}. Finally, we review E.Markman's approach to proving the integrability of the Poisson structure.

      \subsection{Symplectic Structure}
One of the fundamental results in Hitchin's seminal papers \cite{hitchin_sd} and \cite{hitchin_sb} is the discovery that
 $\Higgs_{G,D=0}$ is holomorphic symplectic and carries  the structure of an ACIHS, to be discussed later. 
Recall that,  by definition, a quasi-projective algebraic  variety is  holomorphic symplectic  if its smooth
(regular) locus carries a symplectic   structure, which extends to any desingularisation.
In this subsection we review briefly the construction of the symplectic structure on $\Higgs_{G,D=0}$ and 
then return to the general case $D\neq 0$ in the next subsection.

As we saw  in (\ref{reg_st}),
  $\Higgs_{G,0,c}$  contains  a Zariski open subset  $\Higgs_{G,0,c}^o$ which can be identified with the total space of the
vector bundle $\cF = p_{1 \ast}\left(\ad\left(\pi_c^\ast\scP'\right)\otimes p_2^\ast K_X\right)$
on $\Bun^{rs}_{G,c}$.
Furthermore, one can indentify
$\cF$ with $T^\vee \Bun^{rs}_{G,c}$,  the 
cotangent bundle to the smooth locus of the (coarse) moduli space of semi-stable $G$-bundles of topological type $c$,
and   cotangent bundles carry a canonical symplectic structure.

Pointwise, at the class of a pair $(P,\theta)$,
this identification   can be done  as follows.
By Luna's \'etale slice theorem, $\Bun_{G,c}^{rs}$ is, local-analytically near a  regularly stable bundle $P$, isomorphic to
\[
 H^1(X,\ad P)\sslash \Aut (P) = H^1(X,\ad P)\sslash Z(G)= H^1(X, \ad P)
\]
and  $T_{[P]}\Bun_{G,c}^{rs}= H^1(X,\ad P)$.
Next, the 
stability of  $P$ implies  stability of the Higgs pair  $(P,\theta)$ for any $\theta\in H^0(X, \ad P\otimes K_X)$.
However, 
  a choice of symmetric invariant bilinear form $\textrm{Tr}$ on $\fg$ (e.g., the Killing form)
determines an isomorphism $\ad P= \ad P^\vee$, which, when combined   with
  Serre duality, gives an isomorphism $H^0(X,\ad P\otimes K_X)= H^1(X,\ad P)^\vee$. Hence
the (class of the) pair $(P,\theta)$ determines a point in
$T^\vee\Bun_{G,c}^{rs}$.

The complement $\Higgs_{G,D=0}\backslash T^\vee \Bun^{sm}_{G}$ is non-empty:
there exist stable Higgs pairs with unstable underlying bundle.  A concrete example is furnished 
by the uniformising (or Toda) Higgs bundle, see \cite{hitchin_sd}, Example 1.5, or by any Higgs bundle in
the image of the Hitchin section (\cite{hitchin_teich}).
As shown by these very examples, there are  smooth points in this locus, i.e.,  
\[
\left(\Higgs_{G,D=0}\backslash T^\vee \Bun^{rs}_{G}\right)^{rs}\neq \varnothing,
\]
and we would like to
 extend the symplectic structure to the rest of  $\Higgs_{G,D=0}^{rs}$.

By  a variant of Schlessinger's deformation theory developed in \cite{Biswas-Ramanan}, 
the space of infinitesimal deformations of a Higgs bundle $(P,\theta)$ is
$\HH^1(\scC^\bullet_{(P,\theta)})$,
 where $\scC_{(P,\theta)}^\bullet$ is the Biswas--Ramanan
complex
\begin{equation}\label{br}
\xymatrix@1{\scC_{(P,\theta)}^0= \ad P\ar[r]^-{\ad\theta}\ &\ \ad P\otimes K_X = \scC_{(P,\theta)}^1 }.
\end{equation}
In fact, this is a very special case of  Theorem 2.3, \emph{ibid.}  and   we shall discuss the general case in 
subsection \ref{def_th}.
If $[(P,\theta)]\in \Higgs^{rs}_{G,D=0}$ then $\HH^1(\scC^\bullet_{(P,\theta)})= T_{[P,\theta]}\Higgs_{G,D=0}^{rs}$.

Being a (shifted) cone, the complex $\scC_{(P,\theta)}^\bullet$ is an extension of $\ad P$ by $\ad P\otimes K_X[-1]$ and
the long exact sequence of hypercohomology gives a short exact sequence
\begin{equation}\label{h1}
  \xymatrix@1{(0)\ar[r] & \cok h^0(\ad\theta )\ar[r]& \HH^1\left(\scC_{(P,\theta)}^\bullet\right)\ar[r]^-{\pr}& \ker h^1(\ad\theta)\ar[r]& (0)}    .
\end{equation}
Here $h^i(\ad\theta): H^i(\ad P)\to H^i(\ad P\otimes K_X)$ are the natural maps induced by $\ad\theta$.
If $P$ happens to be  stable, equation (\ref{h1}) reduces to
\begin{equation}\label{h2}
  \xymatrix@1{(0)\ar[r] & H^0(X, \ad P\otimes K_X)\ar[r]& \HH^1\left(\scC_{(P,\theta)}^\bullet\right)\ar[r]^-{\pr}& H^1(X, \ad P) \ar[r]& (0)  }.
\end{equation}
Next,   the combination of  $\textrm{Tr}$ and cup product pairing 
\begin{equation}\label{cup_1}
 \HH^1\left(\scC_{(P,\theta)}^\bullet\right)\otimes \HH^1\left(\scC_{(P,\theta)}^\bullet\right)\to  H^1(K_X)\simeq \CC
\end{equation}
induces a skew-symmetric bilinear form $\omega_{(P,\theta)}\in \Lambda^2\left(\HH^1(\scC_{(P,\theta)}^\bullet)^\vee\right)$. 
As shown in \cite{Biswas-Ramanan}, Theorem 4.3, 
this pairing gives rise to
a   symplectic form on $\Higgs_{G,D=0}^{rs}$ which coincides with   the canonical symplectic form
$\omega^{can}=-d\lambda$ on $\tot T^\vee \Bun^{sm}_{G}$. 
In terms of the deformation complex, the Liouville 1-form 
$\lambda_{(P,\theta)}\in \HH^1\left(\scC_{(P,\theta)}^\bullet\right)^\vee$ is given by
\[
 \lambda_{(P,\theta)}(v)=\textrm{Tr }\pr (v)\cap \theta .
\]

The symplectic form $\omega$ determines (and is determined by) a map 
\begin{equation}\label{symp_contr}
 \omega\lrcorner:\ T_{\Higgs^{rs}_{G,D=0}}\longrightarrow T^\vee_{\Higgs^{rs}_{G,D=0}}.
\end{equation}
At $(P,\theta)$ this corresponds to the linear map
$\HH^1\left(\scC^\bullet_{(P,\theta)}\right)\to \HH^1\left(\scC^\bullet_{(P,\theta)}\right)^\vee$
determined by  Grothendieck--Serre duality for hypercohomology. We are going to elaborate on this 
in the next subsection (see (\ref{duality}), (\ref{duality_2})) where we address
 the case $D\neq 0$.

Finally, we recall how to express the symplectic form  in Dolbeault terms. For that, one considers the
(global sections of the total complex of the) Dolbeault resolution of the complex $\scC^\bullet_{(P,\theta)}$.
The deformation theory of a stable Higgs pair is formal (\cite{hbls}, Lemma 2.2). The  Hermite--Yang--Mills
metric on $(P,\theta)$ provides an embedding $\HH^1\left(\scC^\bullet_{(P,\theta)}\right)\subset A^1(\ad P)$, whose
 image consists of harmonic representatives of hypercohomology.
In terms of the type decomposition
\[
  A^1(\ad P)\simeq A^{0}(\ad P\otimes K_X)\oplus A^{0,1}(\ad P),
\]
the symplectic form is given by restricting the pairing
\[
 \omega \left((\eta',\eta''),(\xi',\xi'')\right)=\int_X \textrm{Tr }\left(\eta'\wedge \xi''-\xi'\wedge\eta''\right) =\int_X\textrm{Tr }(\eta'+\eta'')\wedge (\xi'+\xi'')
\]
to the harmonic representatives of $\HH^1(\scC^\bullet)$.

This brings us back to Hitchin's original motivation: if $K\subset G$ is a maximal compact subgroup, and $Q\subset P$
a $K$-reduction, then the space of holomorphic structures on  $Q\times^K G$ is an affine space
modelled on 
$A^{0,1}(\ad Q_{\CC})$.  This torsor is canonically trivialised and identified with $A^{0,1}(\ad P)$ by the holomorphic structure of $P$.
Then
 $A^1(\ad P)$ can be thought of as the total space of the
(weak) cotangent bundle to the space of holomorphic structures on $Q\times^K G$, and Hitchin's original construction of
$\Higgs_{G,D=0}$ was a kind of  infinite-dimensional hyperkaehler
Marsden--Weinstein reduction of the latter. The holomorphic symplectic form is thus a reduction of the canonical (weak) symplectic form
on the product of a vector space with its (weak) dual.

	      \subsection{Deformation theory}\label{def_th}
		  \subsubsection{One-parameter analytic deformations}
The symplectic structure on the moduli space of $K_X$-valued $G$-Higgs bundles  was defined  in terms of the  complex
(\ref{br})
which controls the  infinitesimal deformations of the $K_X$-valued Higgs pair. 
For $K_X(D)$-valued Higgs bundles a similar complex exists. In fact,
Biswas and Ramanan (\cite{Biswas-Ramanan}) have given a  uniform description of the deformation theory of
Higgs bundles with coefficients in an arbitrary vector bundle.
 Before turning to  their  theorem, which is concerned with infinitesimal
deformations, we discuss the global case in the analytic category.

Let  
$\rho: G\to \Aut F$ be a linear representation of an algebraic group $G$, $P$ a principal $G$-bundle, $\rho P = P\times ^G F$ the corresponding associated vector 
bundle and $\phi\in \Gamma(X, \rho P)$.
  Consider now the analytifications of these objects. Let $\Delta=\left\{\epsilon: |\epsilon|<1 \right\}\subset \CC$ be the unit disk, 
 and let $\scP\to X\times \Delta$ be a deformation of $P$, i.e., a  holomorphic $G$-bundle together
with an isomorphism $\left. \scP\right|_{X\times\{0\}}=P$. 
Let the section $\Phi \in \Gamma(X\times \Delta, \scP\times^G F)$ be a 
deformation of $\phi$, i.e, $\left. \Phi\right|_{X\times\{0\}}=\phi$ under the above isomorphism of bundles.
The section $\Phi$ corresponds to a holomorphic map $\sigma_\Phi: \scP\to F$, which is $G$-equivariant, i.e., 
$R_g^\ast \sigma_\Phi = \rho(g^{-1})\circ \sigma_\Phi$. Similarly, $\phi$ corresponds to a $G$-equivariant map
$\sigma_\phi: P\to F$, and $\sigma_\phi= \left. \sigma_\Phi\right|_{P}$.

 Fix an ``admissible'' cover $\fU= \{U_i\}$ of $X$, i.e., one for which  the family $\scP$ is
trivial over $U_i\times \Delta$. Let us also fix trivialisations, i.e., $G$-bundle isomorphisms
$\Psi_i:  \scP_{U_i\times \Delta}\simeq P_{U_i}\times \Delta = p_1^\ast P_{U_i}$, such that
$\left. \Psi_i\right|_{\epsilon=0}=id$.
Then the composition $\Psi_{ij}=\Psi_i\circ \Psi_j^{-1}\in \Aut (P_{U_{ij}}\times \Delta)$ corresponds to a $G$-equivariant
map $\psi_{ij}: P_{U_{ij}}\times \Delta\to G$, defined by
\[
\Psi_{ij}(p,\epsilon)= (p,\epsilon)\cdot \psi_{ij}(p,\epsilon) = (p\cdot \psi_{ij}(p,\epsilon),\epsilon).
\]
The $G$-equivariance of $\psi_{ij}$ is with respect to the conjugation action, i.e., $R_g^\ast\psi_{ij}= Ad(g^{-1})\circ \psi_{ij}$, or, 
pointwise, $\psi_{ij}(p\cdot g,\epsilon)= g^{-1}\psi_{ij}(p,\epsilon) g$.

Taking into account that $\psi_{ij}(p,0)=e\in G$, we see that the derivative
\[
 \dot{\psi}_{ij}= \left. \frac{d}{d\epsilon}\right|_{\epsilon=0}\psi_{ij}= \left. d\psi_{ij}\left( \frac{d}{d\epsilon}\right)\right|_{\epsilon=0}: P_{U_{ij}}\to \fg
\]
satisfies $R_g^\ast \dot{\psi}_{ij}= \ad (g^{-1})\circ\dot{\psi}_{ij}$, and hence determines a section $s_{ij}\in \Gamma(U_{ij},\ad P_{U_{ij}})$. 
The cocycle condition $\Psi_{ij}\circ \Psi_{jk}\circ \Psi_{ki}= id\in \Aut (P_{U_{ijk}})$ translates to
\[
\psi_{ij}\psi_{jk}\psi_{ki}=e : P_{U_{ijk}}\times \Delta\to G,
 \]
 which, in turn, gives
$s_{ij}+ s_{jk}+ s_{ki}=0$, i.e., $\underline{s}= (s_{ij})\in \check{Z}^1_{\fU}(\ad P)$.

It is then easy to see that the holomorphic maps
\[
 \tau_i = \sigma_{\Phi}\circ \Psi_i^{-1}: P_{U_i}\times \Delta\to F
\]
are $G$-equivariant and satisfy $\tau_j =\tau_i\circ \Psi_{ij}$, or, equivalently, $\tau_j = \rho(\psi_{ji})\circ\tau_i$.
Differentiating this condition at $\epsilon=0$ gives that 
\[
t_i = \left. \frac{d}{d\epsilon}\right|_{\epsilon=0}\tau_i: P_{U_i}\to F
\]
satisfy $t_j- t_i = \rho(\dot{\psi}_{ji})(\phi)$, i.e. $\delta^0\underline{t}= \rho(\underline{s})(\phi)$,
where $\delta^0$ denotes the \v{C}ech differential.  

The infinitesimal data associated to the deformation $(\scP,\Phi)$ is encoded in 
the pair $(\underline{s},\underline{t})$, and it is not hard to trace how this data changes when we pass
to an equivalent deformation. 
		  \subsubsection{The deformation functor}
We turn now to infinitesimal deformations. Let $G$, $P$ and $\rho$ be as before, and let $V$ be a vector bundle on $X$.
Denote by $\scF_{(P,\phi)}:\Art\to \sets$  the formal deformation functor of
the  pair  $(P,\phi)$,  where now $\phi\in \Gamma(X, \rho P\otimes V)$. This is a functor from the category of Artin local $\CC$-algebras to the category of sets,
 which assigns to an algebra $A$
the set of iso-classes of deformations of $(P,\theta)$, parametrised by $X\times \spec A$.
In particular, $\scF_{(P,\phi)}(\CC[\eps]/\eps^2)$ is the space of infinitesimal deformation of the pair $(P,\phi)$.

      \begin{thm}[Theorem 2.3, \cite{Biswas-Ramanan}]\label{br_thm}
  There exists a canonical bijection
\[
\scF_{(P,\phi)}(\CC[\eps]/\eps^2)= \HH^1(\scK^\bullet_{(P,\phi)}),
\]
 where $\scK_{(P,\phi)}^\bullet$
is the complex
\begin{equation}\label{br_gen}
 \xymatrix@1{ \scK^0_{(P,\phi)}= \ad P\ar[r]^-{e(\phi)}& \rho P\otimes V = \scK^1_{(P,\phi)}},
\end{equation}
and $e(\phi)(s)= \rho(s)(\phi)$.     
      \end{thm}
  For simplicity, we have not included $\rho$ or $V$ in the notation of the complex.

  \emph{Sketch of proof:}
The theorem is proved by a direct infinitesimal calculation:  if
$\{U_i=\spec A_i\}_i$ is an affine cover of $X$, then  $U_i[\eps]= \spec \left(A_i\otimes \CC[\eps]/\eps^2\right)$
is an affine cover of $X[\eps]=X\times \spec \CC[\eps]/\eps^2$, and we can replace  $\scK^\bullet_{(P,\phi)}$
by  its \v{C}ech resolution. Then elements of $\HH^1(\scK^\bullet_{(P,\phi)})$ are  identified with equivalence classes  
of pairs
\[
(\underline{s},\underline{t}) = \left((s_{ij})_{ij}, (t_i)_i\right) \ \in \oplus_{ij\ } \ad P (U_{ij}) \bigoplus \oplus _i \left(\rho P\otimes V\right)(U_i),
\] 
 which on  double overlaps
satisfy the two conditions
    \begin{enumerate}
     \item $\delta^1 \underline{s}=0$
     \item $e (\phi)(s_{ij})= (\delta^0 \underline{t})_{ij}$. 
    \end{enumerate}
The restrictions $P_i=\left. P\right|_{U_i}$ of  $P$ determine  trivial families $\scP_i=p_{U_i}^\ast P_i$ on $U_i[\eps]$.
The first condition states that the automorphisms $(1+ s_{ij}\eps)\in\textrm{Aut}(\scP_{ij})$ glue these trivial families into a $G$-bundle
$\scP$
on $X[\epsilon]$. The second condition guarantees that the local sections $(\phi +t_i\eps)_i$ glue into  a section of $\rho\scP\otimes V$.
This determines the map $\HH^1(\scK^\bullet_{(P,\phi)})\to \scF(\CC[\eps]/\eps^2)$. 
The map in the opposite direction is obtained by observing that any deformation 
of $P_i$ over $U_i[\eps]$ is trivial. Hence, given a deformation $(\scP,\Phi)$ of $(P,\phi)$ parametrised by $\spec \CC[\eps]/\eps^2$,
we can  fix trivialisations and obtain the corresponding  gluing data $\underline{s}=(s_{ij})$ and $\underline{t}=(t_i)$, where
and $t_i\eps = \Phi_i - \left. p_X^\ast\phi\right|_{U_i[\eps]}$.
\qed

 Eventually, we are  interested in applying this theorem in the case  
$V=K_X(D)$ and
$\rho(s)=-\ad s$. Then the complex $\scK^\bullet_{(P,\theta)}$ takes the form
\begin{equation}\label{br_mero}
 \xymatrix@1{  \ad P\ar[r]^-{\ad\theta }\ &\ \ad P\otimes K_X(D)  },
\end{equation}
and fits in the extension
\begin{equation}\label{cone_mero}
 \xymatrix@1{0\ar[r]& \ad P\otimes K_X(D)[-1]\ \ar[r]& \scK^\bullet_{(P,\theta)}\ \ar[r] & \ad P\ \ar[r]& 0}.
\end{equation}
We can use this exact sequence to calculate the dimension (\ref{higgs_dim}) of $\Higgs_{G,D}$.
Indeed, taking Euler characteristics gives
\[
 -\chi \left(\ad P\otimes K_X(D)\right) - \chi \left(\scK^\bullet_{(P,\theta)}\right) + \chi \left(\ad P\right) =0,
\]
and hence, by Hirzebruch--Riemann--Roch, 
\begin{equation}\label{dim}
 \dim \Higgs_{G,D}= - \chi\left(\scK^\bullet_{(P,\theta)}\right)= \dim G\deg K_X(D),
\end{equation}
which is also nothing but $ 2\dim \Bun_G + \dim G\deg D$.

By Luna's \'etale slice theorem one can identify an \'etale (or analytic) neighbourhood of $(P,\theta)$
as
\[
 \HH^1 (\scK^\bullet_{(P,\theta)})\sslash \Aut (P,\theta),
\]
which for regularly stable pairs reduces to
 \[
\HH^1(\scK^\bullet_{(P,\theta)})\simeq T_{(P,\theta)}\Higgs^{rs}_{G,D,c}.
\]
	\subsection{Digression on duality}\label{duality_sect}
To handle  the Poisson structure on $\Higgs_{G,D,c}$  we need a small amount  of duality theory (which was already used implicitly in 
(\ref{cup_1})).
 Given a length-$(n+1)$ complex of locally free sheaves 
\[
 \left(F^\bullet, d_\bullet\right) = \left(\xymatrix@1{F^0\ar[r]^{d_0}&  F^1\ar[r]^-{d_1}&\ldots \ar[r]^{d_{n-1}}& F^n }\right),\
F^\bullet = \bigoplus_{k\in\ZZ} F^k\left[-k\right],
\]
 let us denote by $\widehat{F}^\bullet$  its (naive) dual complex, i.e., the $\mhom$-complex
(graded as usual) between the complex $F^\bullet$ and the complex $\cO_X$ (concentrated in degree zero):
\[
\widehat{F}^\bullet =\shom^\bullet_{D_X} (F, \cO_X)= \bigoplus_{k\in \ZZ} \shom_{\cO_X} (F^{-k}, \cO_X)[k].
\] 
This complex is concentrated 
in degrees
$(-n)$ to $0$ and
has differentials which are the duals of the respective differentials of $F^\bullet$:
\[
\xymatrix@1{\delta_{-k}= d_{k-1}^\vee: & \widehat{F}^{- k}= \shom (F^k,\cO_X)\ar[r]& \shom (F^{k-1},\cO_X)= \widehat{F}^{-k+1}}.
\] 
To place $\widehat{F}^\bullet$  in  non-negative degree we
shift it by $n$ positions to the right  and denote the new complex by $\check{F}^\bullet$, i.e.,  $\check{F}^\bullet= \widehat{F}^{\bullet}[-n]$.
Then Grothendieck--Serre duality (Theorem 3.12, \cite{huybrechts_fm}), which in this case is just Serre duality for hypercohomology,
 tells us that for all $i\in\ZZ$
\begin{equation}\label{duality}
\HH^i (F^\bullet)^\vee = \HH^{-i} \left(\widehat{F}^{\bullet}\otimes K_X[1]\right) = \HH^{n+1-i}\left(\check{F}^\bullet\otimes K_X\right).
\end{equation}
The duality can be made  explicit as follows.
The contractions  $F^k\otimes \shom(F^k,\cO_X)\to \cO_X$
give rise to a linear map
\[
 \left(F^\bullet\bigotimes \left(\check{F}^\bullet\otimes K_X\right)\right)_n = \bigoplus_{k=0}^n F^k\otimes \shom(F^{k},\cO_X)\otimes K_X \longrightarrow K_X
\]
and hence to a morphism of complexes 
\begin{equation}\label{duality_tensor}
 F^\bullet\bigotimes \left(\check{F}^\bullet\otimes K_X\right)
\longrightarrow K_X[-n].
\end{equation}
This is indeed a morphism of complexes, since
 the tensor product complex, being the total complex of a double complex, has differential obtained from the tensor product of  the differentials of the
two complexes, with alternating signs.
The morphism (\ref{duality_tensor})  induces 
 a map on cohomology
\[
 \HH^{n+1} \left(F^\bullet\otimes \left(\check{F}^\bullet\otimes K_X\right)\right) \longrightarrow \HH^{n+1}(K_X[-n])= H^1(X,K_X)\simeq \CC.
\]
Then the duality pairing corresponding to (\ref{duality}) is the  composition of this map with the cup product pairing
\[
 \HH^i (F^\bullet)\otimes \HH^{n+1-i}\left(\check{F}^\bullet\otimes K_X\right) \longrightarrow \HH^{n+1} \left(F^\bullet\otimes \left(\check{F}^\bullet\otimes K_X\right)\right).
\]

The case of interest for us is
 the complex (\ref{br_mero}) which has length two ($n=2$)
and hence
\begin{equation}\label{duality_2}
 \HH^i(\scK^\bullet_{(P,\theta)})^\vee = \HH^{2-i} (\check{\scK}^\bullet_{(P,\theta)}\otimes K_X)
\end{equation}
where
\[
 \xymatrix@1{\check{\scK}^\bullet_{(P,\theta)}\otimes K_X:\ &  \ad P^\vee\otimes \cO_X(-D)\ar[r]^-{(\ad \theta)^\vee}\ &\ \ad P^\vee\otimes K_X    }.
\]
Using the chosen Ad-invariant symmetric bilinear form on $\fg$, we identify $\check{\scK}^\bullet_{(P,\theta)}\otimes K_X$
with the complex
\[
 \xymatrix@1{  \ad P\otimes \cO_X(-D)\ar[r]^-{- \ad \theta }\ &\ \ad P\otimes K_X     }.
\]

  \begin{rmk}
Naturally, changing the sign of the differential in the complex (\ref{br_gen}) gives an isomorphic complex, but we, nonetheless, make some comments on 
the sign choices involved.   
We use the (standard)
 convention for the $0$-th \v{C}ech differential $(\delta^0\underline{t})_{ij}= t_j- t_i$ and the (fairly standard) differential 
$e(\phi) + (-1)^{i+1}\delta^j: \check{C}_{ij}\to \check{C}_{i+1,j}\oplus \check{C}_{i,j+1}$ for the \v{C}ech double complex.
This forces 
our definition of $e(\phi)$
 in Theorem \ref{br_thm} to differ by sign from the one in \cite{Biswas-Ramanan}.
On the other hand, in equation (\ref{br_mero}) we have used a sign, corresponding to $-\ad: \fg\to \mend(\fg)$, for the following reason.
In view of  duality, it would have been more natural if Higgs fields were defined using 
the co-adjoint representation, i.e., as sections of $(\ad P)^\vee\otimes K_X$, as in \cite{Biswas-Ramanan}, \S 4. In that case
the deformation complex would satisfy $\check{\scK}^\bullet_{(P,\theta)}\otimes K_X= \scK^\bullet_{(P,\theta)}$,
without the need to choose an invariant symmetric bilinear form. Since we stick, however, to  the standard definition of a Higgs field,
and $\textrm{Tr}$ identifies $(\ad\theta)^\vee$ with $-\ad\theta$, the above sign choice is forced unto us.
  \end{rmk}

  \begin{Prop}
Let $(P,\theta)$ be a $K_X(D)$-valued $G$-Higgs bundle and
 $\fU=\{(U_i)\}$  an acyclic cover of $X$. 
Let $\alpha = \left[(\underline{s},\underline{t}) \right]$ and
$\beta=\left[(\underline{\sigma},\underline{\tau}) \right]$ be hypercohomology classes in
$ \HH^1(\scK_{(P,\theta)}^\bullet)$ and $\HH^1(\check{\scK}^\bullet_{(P,\theta)}\otimes K_X)$, respectively, 
with \v{C}ech representatives
\[
(\underline{s},\underline{t})\in \check{C}^1_{\fU}(\ad P)\oplus \check{C}_{\fU}^0(\ad P\otimes K_X(D))
\]
and
\[
(\underline{\sigma},\underline{\tau})\in  \check{C}^1_\fU(\ad P (-D))\oplus \check{C}^0_\fU(\ad P\otimes K_X).
\]
Then the
duality pairing
\[
 \HH^1(\scK_{(P,\theta)}^\bullet)\otimes \HH^1(\check{\scK}^\bullet_{(P,\theta)}\otimes K_X) \longrightarrow H^1(X,K_X)
\] 
corresponding to  (\ref{duality_2}) maps $\alpha\otimes\beta$ to $\left[\underline{c}\right]\in H^1(X,K_X)$,
where $\underline{c}= (c_{ij})$ is given by
\begin{equation}\label{duality_explicit}
c_{ij} = \textrm{Tr}(t_i,\sigma_{ij})-\textrm{Tr}(s_{ij},\tau_j)  \in K_X(U_{ij}). 
\end{equation}
  \end{Prop}
\emph{Proof:} 
In view of (\ref{duality_2}) and (\ref{duality_tensor}), to construct the explicit pairing we need to have an explicit description 
of cup product in hypercohomology. To simplify notation, set $F^\bullet = \scK^\bullet_{(P,\theta)}$ and $G^\bullet= \check{\scK}^\bullet_{(P,\theta)}\otimes K_X$.
Moreover, to avoid confusion, let us be explicit about tensor products of complexes, e.g., write $\textrm{tot}^\bullet (F^\cdot\otimes G^\cdot)$, 
rather than just $F^\bullet\otimes G^\bullet$. Finally, 
 let us write $t\cC^\bullet(F^\cdot)$ for the (total) \v{C}ech complex of $F^\bullet$ with respect to the
cover $\fU$. This is the complex of vector spaces, whose $k$-th term is
\[
t\cC^k(F^\bullet) = \check{C}^k_\fU(F^0)\oplus \check{C}^{k-1}_\fU (F^1) = \oplus_{i_0\ldots i_k}F^0(U_{i_0\ldots i_k})\bigoplus 
\oplus_{i_0\ldots i_{k-1}} F^1(U_{i_0\ldots i_{k-1} })
\]
and whose $k$-th differential is $\begin{pmatrix}
                                   \delta^k &0\\
                                   \ad\theta &-\delta^{k-1}
                                  \end{pmatrix}: t\cC^k(F^\bullet)\to t\cC^{k+1}(F^\bullet)$.
Hypercohomology is computed as
\[
 \HH^n(F^\bullet) = H^n (t\cC^\bullet(F^\cdot)),
\]
and by K\"unneth formula we have
\[
 \HH^n (F^\bullet)\otimes \HH^m(G^\bullet)= H^n (t\cC^\bullet(F^\cdot))\otimes H^m (t\cC^\bullet(G^\cdot)) 
= H^{n+m}\left(\textrm{tot}^\bullet\left( t\cC^\cdot(F^\cdot)\otimes t\cC^\cdot (G^\cdot)\right)\right).
\]
The cup product map 
\[
 \HH^n(F^\bullet)\otimes \HH^m(G^\bullet) \longrightarrow  \HH^{n+m} (\textrm{tot}^\bullet\left( F^\cdot\otimes G^\cdot\right))
\]
is induced by a morphism of complexes of vector spaces
\[	
 \textrm{tot}^\bullet\left(t\cC^\cdot(F^\cdot)\otimes t\cC^\cdot (G^\cdot)\right) \longrightarrow t\cC^\bullet (\textrm{tot}\left(F^\cdot\otimes G^\cdot\right)),
\]
\[
 \alpha\otimes\beta \longmapsto \alpha\cup \beta,
\]
where, if $\deg \alpha=n$, $\deg \beta =m$, one sets
\[
 \left(\alpha\cup\beta\right)_{i_0\ldots i_p}= 
\sum_{r=0}^p (-1)^{r(m-(p-r))}\alpha_{i_0\ldots i_r}\otimes \beta_{i_r\ldots i_p}\in \check{C}^p_\fU(\textrm{tot}^{n+m-p}\left(F^\bullet\otimes G^\bullet\right)).
\]
For a discussion of the sign one can consult Deninger (\cite{deninger_sign}), Deligne (\cite{Deligne_propre}) or 
 A.de Jong's unpublished notes on algebraic de Rham cohomology.

If $n=m=1$ and we consider an element in the image of the K\"unneth map, i.e., 
$\alpha\otimes \beta\in t\cC^1(F^\bullet)\otimes t\cC^1(G^\bullet)$, then
\[
 \alpha\cup\beta \in \check{C}^2_\fU((F^\bullet\otimes G^\bullet)^0)\bigoplus \check{C}^1_\fU((F^\bullet\otimes G^\bullet)^1)
\bigoplus \check{C}^0_\fU((F^\bullet\otimes G^\bullet)^2)
\]
has a component of \v{C}ech degree 1 equal to
\[
 \left(\alpha\cup\beta\right)_{i_0i_1} = \alpha_{i_0}\otimes \beta_{i_0i_1}- \alpha_{i_0i_1}\otimes\beta_{i_1} = t_{i_0}\otimes \sigma_{i_0i_1}
- s_{i_0i_1}\tau_{i_1} \in \check{C}^1_\fU((F^\bullet\otimes G^\bullet)^1)
\]
Projecting onto this component and applying trace gives the claimed formula.
\qed

	  \subsection{Poisson structure}\label{poisson_sect}
		\subsubsection{Generalities}
A Higgs bundle is  a decorated principal bundle: a pair, consisting of a principal bundle $P$ and a section of the  vector bundle
$\rho P\otimes V$.
From this viewpoint, it is only natural to ask whether the  symplectic structure on $\Higgs_{G, D=0}$  persists when one 
varies the representation $\rho$
or the   coefficient bundle $V$.

The symplectic structure   was constructed from two ingredients:
the identification
\[
\scF_{(P,\theta)}(\CC[\eps]/\eps^2)= \HH^1(\scC^\bullet_{(P,\theta)})= T_{(P,\theta)}\Higgs_{G, D=0}^{rs},
\]
due to Theorem \ref{br_thm}, and the
 natural skew pairing (\ref{cup_1}), which gives the isomorphism (\ref{symp_contr}).
While the first ingredient makes sense in   general (after replacing $\scC^\bullet_{(P,\theta)}$ with $\scK^\bullet_{(P,\theta)}$),
the second one 
relies substantially  on the fact that the coefficient bundle is the dualising sheaf $K_X$,
and thus we cannot expect the moduli space to be symplectic for arbitrary $V$.
However, if   $V\simeq K_X(D)$, $D> 0$,
     then it turns out that
an analogue of the dual of (\ref{symp_contr}) still exists. While it may fail to be  everywhere of maximal rank,
it still satisfies the appropriate integrability condition.
More precisely,    Bottacin \cite{bottacin} (for  $G=SL_n(\CC)$ and $G= GL_n(\CC)$) and Markman 
(\cite{markman_thesis}, \cite{markman_sw})  showed that, 
 whenever nonempty, $\Higgs_{G,D,c}$ is a holomorphic Poisson variety.

We recall the definition of Poisson structure below but
refer to \cite{arnold_givental}, \cite{weinstein_local_structure} and
\cite{donagi_markman} for details.

Let $M$ be a smooth analytic (or quasi-projective algebraic) variety and $\Pi \in H^0(M,\Lambda^2 T_M)$ a bivector
(field). It determines a $\CC$-bilinear skew-symmetric pairing $\cO_M\ctimes\cO_M\to \cO_M$
\[
 (f,g)\longmapsto \{f,g\}:= \left(df\wedge dg\right)\lrcorner \Pi,
\]
which is a $\CC$-derivation in each entry. Hence
to a (local) function $f\in \cO_M(U)$, $U\subset M$, we can associate a 
$\CC$-derivation of the $\CC$-algebra $\cO_M(U)$, called
its \emph{Hamiltonian vector field} 
\[
X_f = \left\{f,\ \right\} = df \lrcorner \Pi\ \in T_M(U)=\underline{Der}_{\CC}(\cO_M)(U).
\]
We say that $\Pi$ is a \emph{Poisson structure} if this pairing  endows $\cO_M$ with the structure of a sheaf
of Lie algebras. The bracket  $\{f,g\}$ is then called \emph{the Poisson bracket} of $f$ and $g$. The variety 
 $M$ is said to be \emph{Poisson} if it admits a Poisson structure.

The Jacobi identity is equivalent to the requirement that the ``adjoint representation'' $f\mapsto \{f,\ \}$
be not only $\CC$-linear, but also a  homomorphism $\cO_M\to \underline{Der}_\CC(\cO_M)$ of sheaves of Lie algebras, i.e., 
$[X_f,X_g]= X_{\{f,g\}}$.
This 
 can also be
phrased as the vanishing of the Schouten bracket of $\Pi$ with itself. 

On a Poisson variety $(M,\Pi)$ one has an obvious
sheaf homomorphism $\Psi: T_M^\vee\to T_M$, namely,   $\alpha \mapsto \alpha\lrcorner \Pi$.   
It  gives rise to  a stratification of
$M$ by submanifolds $M_k$, $k$-even, such that $\rk \Psi_{M_k} =k$. Then  $\left(M_k, \left. \Pi\right|_{M_k}\right)$
is Poisson.
 The strata are further foliated (\cite{weinstein_local_structure}), local-analytically,
by the $k$-dimensional integral leaves $S$ of the distribution $\left. \Psi\right|_{M_k} (T^\vee_{M_k})\subset T_{M_k}$. 
If  $S\subset M_k$ is a symplectic leaf, then
 $\left. \Pi\right|_{S}$ is   a symplectic structure on it. The leaves can also be identified as the level sets of the Casimir functions, i.e., the functions
$f\in H^0(M,\cO_M)$ with $X_f=0$.

One of the best-known examples of Poisson structure is the \emph{Kostant--Kirillov} Poisson structure. If 
$\cG$ is a Lie group with Lie algebra $\bg$, then its (linear) dual $\bg^\vee$ carries a Poisson bracket
\[
 \{f,g\}_\alpha =\alpha([df_\alpha, dg_\alpha]). 
\]
The group $\cG$ acts on $\bg^\vee$ via the coadjoint representation, and the coadjoint orbits are the symplectic
leaves of the Kostant--Kirillov Poisson structure. The Casimir functions are the invariants
$\CC[\bg]^\cG\subset H^0(\bg^\vee,\cO_{\bg^\vee})$.

	      \subsubsection{A bivector on $\Higgs_{G,D}$}
The moduli space $\Higgs_{G,D,c}$ carries a canonical bivector $\Pi$ which can be described entirely in terms of
homological algebra and is a natural candidate for a Poisson structure. We discuss it below, roughly along the lines of
\cite{bottacin} \S 3, \cite{markman_thesis} \S 7.2 and
\cite{donagi_markman}, \S 5.4. Since these references deal exclusively with the case of $G=GL_n(\CC)$, we spell certain
points in more detail, but see also \cite{markman_sw}.

The canonical inclusion $s: \cO_X(-D) \hookr \cO_X$ 
induces   a  morphism of complexes
\begin{equation}\label{incl}
\xymatrix@1{I_s:\ & \check{\scK}^\bullet_{(P,\theta)}\otimes K_X \ar[r]&\ \scK^\bullet_{(P,\theta)}(-D)\ar[r]\ &\ \scK^\bullet_{(P,\theta)} }.
\end{equation}
Assuming that $(P,\theta)$ is regularly stable, the induced map on $\HH^1$ gives
\begin{equation}\label{poisson_1}
\xymatrix@1{ \Psi_{(P,\theta)}= \HH^1 (I_s):\ & T_{(P,\theta)}^\vee\Higgs^{rs}_{G,D,c}  \ar[r]&\
T_{(P,\theta)}\Higgs^{rs}_{G,D,c}= \HH^1(\scK^\bullet_{(P,\theta)})},
\end{equation}
which is easily seen to be skew-adjoint. This map corresponds to an element
\begin{equation}\label{poisson_2}
    \Pi_{(P,\theta)}\in \Lambda^2 \left(\HH^1(\scK_{(P,\theta)}^\bullet)\right) \subset \HH^1(\scK_{(P,\theta)}^\bullet)^{\otimes 2},
\end{equation}
which is our candidate for a Poisson bivector.

Let us give an explicit global description of the Poisson structure over the locus $\Higgs^o_{G,D,c}\subset \Higgs^{rs}_{G,D,c}$
 (see \ref{reg_st}), assuming $g\geq 2$. This locus is a vector bundle
\[
 \xymatrix@1{\tot\cF= \Higgs^o_{G,D,c}\ar[r]^-{p_\cF} & \Bun_{G,c}^{rs}},
\]
and  so  the isomorphism $p_\cF^\ast \cF^\vee = \Omega^1_{\Higgs^o/\Bun}$ gives
\begin{equation}\label{rel_tan_1}
 \xymatrix@1{(0)\ar[r] & p_\cF^\ast\cF\ar[r] & T_{\Higgs^o_{G,D,c}}\ar[r]^-{d p_\cF}& p_\cF^\ast T_{\Bun^{rs}_{G,c}} \ar[r]& (0)  }
\end{equation}
for the relative tangent sequence on $\tot \cF$.
Since $\cF = p_{1 \ast}\left(\ad\left(\pi_c^\ast\scP'\right)\otimes p_2^\ast L\right)$, 
upon restriction to a point $(P,\theta)\in \tot \cF$ this sequence becomes
\begin{equation}
 \xymatrix@1{(0)\ar[r]& H^0(X,\ad P\otimes L)\ar[r]&\ \HH^1(\scK^\bullet_{(P,\theta)})\ar[r]&\ H^1(X,\ad P)\ar[r]& (0)   }.
\end{equation}
This  is nothing but the degree-1 piece of the cohomology sequence of (\ref{cone_mero}). If $D=0$,  this is the 
sequence (\ref{h2}).

 Consider again the diagram (\ref{universal}) and recall that we have an
adjoint Poincar\'e family of $G^{ad}$-Higgs bundles
$(\scQ,\Theta)$ on $\Higgs^o_{G,D,c}\times X$. Here $\Theta\in H^0(\tot \cF, p_\cF^\ast\cF)$ is the tautological section
and
 $\scQ = (\pi_c'\circ p_\cF')^\ast \scP'$. We then have at our disposal the universal Biswas--Ramanan complex on 
$\Higgs^o_{G,D,c}\times X$
\begin{equation}\label{br_univ}
 \xymatrix@1{\scK_{(\scQ,\Theta)}^\bullet:\ & \ad\scQ\ar[r]^-{\ad\Theta} \ &\ \ad\scQ\otimes p_X^\ast L   }
\end{equation}
with its cone sequence
\begin{equation}\label{univ_cone}
  \xymatrix@1{(0)\ar[r]& \ad\scQ\otimes p_X^\ast L \left[-1\right]\ar[r]& \scK^\bullet_{(\scQ,\Theta)} \ar[r]& \ad\scQ\ar[r]& (0)  }.
\end{equation}
The first (hyper-)derived image of $p_{1\ast}$ applied to (\ref{br_univ}) gives the tangent sheaf
\[
T_{\Higgs_{G,D,c}^o}= \bR^1 p_{1\ast} \scK^\bullet _{(\scQ,\Theta)},
\]
 while  applying $p_{1\ast}$ to
(\ref{univ_cone}) gives the relative tangent sequence (\ref{rel_tan_1}).
Similarly,
\[
T^\vee_{\Higgs^o_{G,D,c}} = \bR^1 p_{1\ast}\left(\check{\scK}^\bullet_{(\scQ,\Theta)}\otimes p_2^\ast K_X\right).
\]
Then we have the relative analogue of the map (\ref{incl}), 
\[
 \xymatrix@1{ I_s:\ &\ \check{\scK}^\bullet_{(\scQ,\Theta)}\otimes p_2^\ast K_X\ar[r]&\  \scK^\bullet_{(\scQ,\Theta)}},
\]
and
\begin{equation}\label{poisson_3}
 \xymatrix@1{ \Psi =  \bR^1 p_{1\ast}(I_s):&  T_{\Higgs^o_{G,D,c}}^\vee\ar[r]&\   T_{\Higgs_{G,D,c}^o}}
\end{equation}
determines a bivector $\Pi\in \Lambda^2\left(T_{\Higgs^o_{G,D,c}}\right)$ which restricts to
(\ref{poisson_2}) at each pair $(P,\theta)$ for which $P$ is  regularly stable.
      \subsubsection{Integrability of the Poisson bivector}
The bivector $\Psi$  does indeed determine a Poisson structure, but this is not easy to prove.
In \cite{bottacin} \S 4.6 (see also \S 4.2) 
this was done by a direct cocycle calculation
 rooted in the fact that
the total space of the dual of a Lie algebroid carries a canonical Poisson structure.

 Markman employed in (\cite{markman_thesis}, \cite{markman_sw})  
a different strategy.
He started by considering  $G$-bundles with framing  along the divisor $D$.
The group of framings (level group) acts on this space and the action lifts to its cotangent bundle.
Over a certain Zariski open subset of the latter  the action is free and the quotient can be identified with
an open subset of $\Higgs_{G,D,c}$. By general properties of Marsden--Weinstein reduction, the canonical symplectic structure
on the cotangent bundle to the moduli space of framed bundles 
descends to a Poisson structure on the reduced space. Markman then  verified that this Poisson structure coincides
with the one induced by the general hypercohomological argument above. Consequently, $\{\Pi,\Pi\}=0$ everywhere.
We review   Markman's construction in the next subsection.

	  \subsection{Framed Bundles and Markman's construction} 

		\subsubsection{Jet Schemes}
If $D\subset X$ is a (possibly non-reduced) divisor and $P\to X$ is a principal $G$-bundle, then 
data of framing of $P$ along $D$
is encoded in points of certain jet schemes of $G$.
In order to make the exposition self-contained we 
 recall here the definition and the
 explicit description of jet schemes of affine varieties.

To any given    scheme  $\cY$ (of finite type over $\CC$)
one  can associate, for any $n\in \NN$,  the functor
\[
 \xymatrix@1{\mhom_{Sch_\CC}(\ \underline{ \hspace{5pt}} \times \spec \CC[t]/t^{n+1}, \cY):& \textrm{Sch}_{\CC}^{op}\ar[r]& \sets }.
\]
 This functor is  representable (\cite{ein_mustata_jets}, \S 2) by a
 scheme $\cY_n$ (of finite type, over $\CC$),
 known  as \emph{the $n$-th jet scheme of $\cY$}. It is easy to see that
 $\cY_0= \cY$, that $\cY_1=\tot T_\cY=\underline{\spec} \Omega^1_\cY$ and that
 there are natural maps $\cY_n\to \cY_{n-1}$.
Notice that by
 definition, for any  $\CC$-algebra $A$ one has
\[
 \mhom_{Sch_\CC} (\spec A[t]/t^{n+1}, \cY) = \mhom_{Sch_\CC} (\spec A, \cY_n),
\]
and in particular,  $\cY_n(\CC)$ is  identified with
$\mhom_{Sch_\CC}(\spec \CC[t]/t^{n+1},\cY)$,   the set of $n$-jets of paths into $\cY$.

We describe now $\cY_n$   for an embedded affine variety 
 $\cY=\spec R\subset \mathbb{A}^N$, where  $R= \CC[x_1,\ldots, x_N]/\fA$.
 For that we shall exhibit  a $\CC$-algebra $R_n$ together with an  isomorphism
\[
 \mhom_{alg} (R, A[t]/t^{n+1}) = \mhom_{alg} (R_n, A),
\]
 functorial in $A$, and set   $\cY_n=\spec R_n$.
We consider first  $\cY=\mathbb{A}^N_\CC$, and claim that
$\mathbb{A}^N_n= \mathbb{A}^{N(n+1)}$. Indeed,  $\CC$-algebra homomorphisms
 $  \spec\CC [x_1,\ldots, x_N]\to A$ are in bijection with  matrices
$\textrm{Mat}_{N\times (n+1)}(A)$, since a homomorphism $\phi$ is specified by 
$\phi(x_i)= \sum_{k=0}^n M_{ik}t^k$, $M_{ik}\in A$.
But  such matrices are also in bijection with  algebra homomorphisms $\CC[y_{10},\ldots,y_{1n},\ldots, y_{N0},\ldots, y_{Nn}]\to A$
via $M\mapsto \psi$,   $\psi(y_{ik})=  M_{ik}$.
Suppose next  that $\fA = (f_1, \ldots, f_p)\neq (0)$. The homomorphisms
$R  \to A[t]/t^{n+1}$ are precisely those   homomorphisms 
 $\phi: \CC[x_1,\ldots, x_N]\to A[t]/t^{n+1}$,  which
factor through the quotient, i.e., 
$f_l(\phi(x_1),\ldots, \phi(x_N))=0\in A[t]/t^{n+1}$, $1\leq l\leq p$.
Expanding the
 latter  gives $\sum_{k=0}^n g_{lk}((M_{ij}))t^k=0$, for some polynomials
$g_{lk}\in \CC[y_{10},\ldots, y_{Nn}]$. We then set
$R_n = \CC[y_{10},\ldots, y_{1n},\ldots, y_{Nn}]/(g_{lk})$
and define $\psi: R_n\to A$ by
 $\psi(y_{ik})= M_{ik}$ as before.
 The jet scheme  $\cY_n\subset \mathbb{A}^{N(n+1)}= \mathbb{A}^{N}_n$
is cut out by the $g_{lp}$.  

For non-affine $\cY$, the jet scheme is constructed by gluing the jet schemes of affine patches.
As far as general properties of jet schemes go, we only mention that 
the  natural maps $\cY_n\to \cY_{n-1}$  are $\mathbb{A}^{\dim \cY}$-bundles, and in particular, the non-singularity
of $\cY=\cY_0$ implies the non-singularity of $\cY_n$ for all $n\in \NN$. Moreover, the assignment
$\cY\mapsto \cY_n$ is functorial in $\cY$ and hence gives rise to an endofunctor of the category of schemes of finite type
(over $\CC$, or any algebraically closed field).

We are interested in jet schemes of (affine) algebraic groups.
 As it is easy to see,  $\left(GL_N\right)_n = GL_N(\CC[t]/t^{n+1})$, and by  the above description,  an affine embedding 
$G\subset GL_N(\CC)$ determines an embedding of the corresponding jet scheme $G_n\subset GL_N(\CC[t]/t^{n+1})$.
		\subsubsection{Framed bundles}
Let $D=\sum _{i=1}^s n_i q_i$ be a (sufficiently positive) divisor on $X$, $\bG$ the $\CC$-scheme, corresponding to our
simple group $G= \bG(\CC)$, and
 let $\widetilde{\bG}_D$ stand for the (group) scheme of maps from $D$ to $\bG$. We recall that it is defined as 
the $\CC$-scheme representing the functor
\[
 \xymatrix@1{ \mhom_{Sch_\CC}(\underline{\hspace{5pt}}\times D, \bG):& Sch^{op}_{\CC}\ar[r]& \sets },
\]
so, for any $\CC$-algebra, $A$, we have
\[
   \mhom_{Sch}(\spec A, \widetilde{\bG}_D) = \mhom_{Sch}(\spec A\times D, \bG) = \mhom_{Alg}(H^0(\cO_\bG), A\otimes H^0(\cO_D)).
\]
We  then see that the group of $\CC$-points is a product of the respective jet schemes of $G$, i.e.,  
\[
\widetilde{G}_D := \widetilde{\bG}_D(\CC)= \mhom_{Alg}(H^0(\cO_\bG), H^0(\cO_D))= \prod_{i=1}^s G_{n_i-1}.
\]
The \emph{level group} is defined to be the quotient $G_D = \widetilde{G}_D/ Z(G)$, where the centre
$Z(G)$ is embedded diagonally.

A \emph{framed bundle}, sometimes also called \emph{a bundle with level-$D$ structure} is a pair $(P,\eta)$,
where $P\to X$ is a principal $G$-bundle and $\eta$ is a trivialisation of $P$ at $D$, i.e., an isomorphism
$\eta: \left. P\right|_D\simeq D\times \bG$ of $G$-bundles.
We define 
\[
\misom((P_1,\eta_1),(P_2,\eta_2))\subset \misom (P_1,P_2)
\]
  as the set of isomorphisms $f:P_1\simeq P_2$, satisfying
$\eta_1=\eta_2\circ f_D$. There is a natural action of $G_D$ on the set of isomorphism classes of framed bundles, 
namely, $g\cdot [(P,\eta)]= [(P,\widetilde{g}\circ \eta)]$, where $\widetilde{g}\in \widetilde{G}_D$ is a lift
of $g\in G_D$, and $\textrm{Stab} [(P,\eta)]= \textrm{Im}(\Aut P\to\Aut P_D)^{op}/Z(G)$.

We denote by $\cP = \cP(G,D,c)$ the smooth locus of the moduli space of isomorphism classes of 
stable framed $G$-bundles of topological type $c$. We make some comments on the r\^ole of $D$ in the definition of
stability, and refer to \cite{markman_thesis}, \cite{markman_sw} and \cite{seshadri_asterisque}
for more details.

Let  $\delta:= \deg D$. A vector bundle $E\to X$ is called \emph{$\delta$-stable}, if for any proper subbundle $F\subset E$, 
one has $\frac{\deg F-\delta}{\rk F}< \frac{\deg E-\delta}{\rk E}$, and one defines similarly  
$\delta$-semistability. A framed \emph{vector} bundle $(E,\eta)$ is (semi)stable, if $E$ is $\delta$-(semi) stable. 
 It is clear that
if $E$ is stable, then it is $\delta$-stable for any $\delta\geq 0$, and if $E$ is semi-stable, it is $\delta$-stable for any
$\delta>0$. 
Seshadri in \cite{seshadri_asterisque} (part 4) constructed a projective coarse moduli space of
semi-stable framed coherent sheaves.
A framed
$G$-bundle $(P,\eta)$ shall be called  (semi-)stable, if $\ad P$ is $\delta$-(semi-)stable.  
  By cocycle calculation it is not hard to see that the tangent space
$T_{\cP, (P,\eta)}= H^1(X, \ad P\otimes \cO(-D))$: these are the infinitesimal deformations of $P$ which preserve the framing, i.e., 
vanish along $D$. Serre duality implies that $\tot T^\vee \cP(\CC)$ consists of (classes of) triples $(P,\eta,\theta)$,
$\theta\in H^0(\ad P\otimes K_X(D))$. Such triples, consisting of a Higgs bundle and a framing of the underlying $G$-bundle
 can be called \emph{framed Higgs bundles}.
	      \subsubsection{Symplectic reduction}
The action of the level group on $\cP$ lifts naturally to the $T^\vee\cP$, and the lifted action is given by
$g\cdot [(P,\eta,\theta)]= [(P,\widetilde{g}\circ\eta, \theta)]$. 
After considering the homomorphism $\Aut P\to \Aut (\ad P)$,
Markman's Lemma 6.7 (\cite{markman_thesis}),
applied to $\ad P$ implies the stabilisers of the $G_D$ action on $\cP$ and $T^\vee\cP$, are, respectively
\[
 \textrm{Stab}([P,\eta]) = \Aut (P)^{op}/Z(G)\subset G_D
\]
and
\[
 \textrm{Stab}([P,\eta,\theta])= \Aut (P,\theta)^{op}/Z(G)\subset G_D.
\]
Correspondingly, the action of $G_D$ on the
 loci $\cP^o =\left\{[(P,\eta)], [P]\in \Bun^{rs}_{G,c}\right\}\subset\cP$ and
$(T^\vee\cP)^o= \left\{[(P,\eta,\theta)], [(P,\theta)]\in \Higgs^{rs}_{G,D,c}\right\}\subset T^\vee\cP$ is free, and these are
principal $G_D$-bundles over $\Bun_{G,c}^{rs}$ and $\Higgs_{G,D,c}^{rs}$, respectively.

Lifted actions on cotangent bundles give rise to a very special geometry, as we now recall, 
following \cite{arnold_givental} and \cite{donagi_markman}. Not only is the manifold
 $M:= \tot T^\vee\cP$  symplectic, but, moreover, 
 the action $G_D\times M\to M$ is Poisson, and  there exists 
  a canonical $G_D$-equivariant
moment map $\mu: M\to \fg_D^\vee:= \textrm{Lie}G_D^\vee$. We explain these properties
briefly.

Let
$\ba: \fg_D\to H^0(M,T_M)$ denote the infinitesimal action map, assigning to each $\xi\in \fg_D$
the corresponding ``fundamental vector field''. The action of $G_D$ on $M$ 
is  \emph{Hamiltonian} if $\textrm{Im}(\ba)$ consists of Hamiltonian vector fields:
for any $\xi\in \fg_D$, there is a global function $f\in H^0(\cO_M)$, such that
 $\ba(\xi)= X_f$.
The action is \emph{Poisson}, if it is Hamiltonian and if the hamiltonian functions for the different $\xi\in \fg_D$
can be chosen compatibly, i.e., if there exists a Lie algebra homomorphism
$H: \fg_D\to H^0(M,\cO_M)$ and  $\ba$ factors through it, giving   $\ba(\xi)= X_{H(\xi)}$.

Dually,  $H$ can be thought of as a  \emph{moment map}, i.e., a morphism $\mu: M\to \fg^\vee$.
It is Poisson and  $G_D$-equivariant. For the case that we consider -- a lifted action
of a (connected) group on a cotangent bundle of a manifold -- there is a canonical moment map, \cite{arnold_givental}. Namely, 
if $(u,\theta)\in T^\vee_{\cP,u}$, $\mu(u,\theta) = (d\rho_u)_e^\vee(\theta)$, where $\rho_u: G_D\to \cP$ is the orbit map, 
$\rho_u(g)=g\cdot u$. 

We can identify $\fg_D^\vee=\widetilde{\fg}^\vee_D$ with $\fg^\vee\otimes H^0(\left. K_X(D)\right|_D)$ via  
\[
 \xymatrix@1{H^0(\cO_D)\otimes H^0(\left. K_X(D)\right|_D)\ar[r]&  H^0(\left. K_X(D)\right|_D)\ar[r]^-{\textrm{Res}} & H^1(K_X)\simeq \CC},
\]
where $\textrm{Res}$ is the (first) connecting homomorphism of the long exact cohomology sequence, associated to 
\[
 \xymatrix@1{(0)\ar[r]& K_X \ar[r]&K_X(D)\ar[r] & \left. K_X(D)\right|_{D}\ar[r]& (0) }.
\]
 Then by (\cite{markman_thesis}, Proposition 6.12, \cite{markman_sw}), the moment map is explicated as
\[
\mu([(P,\eta,\theta)])(A)= \textrm{Res } A\left( \eta\circ \left. \theta\right|_{D}\circ\eta^{-1}\right).
\]

If the map $\mu$ were submersive and the quotient $M/G_D$ were to exist, it would carry a canonical Marsden--Weinstein
Poisson structure, whose symplectic leaf through $m\in M$ would be
\[
 \mu^{-1}\left(\bO_{\mu(m)}\right)/ G_D \simeq \mu^{-1}(\mu(m))/ \textrm{Stab}(\mu(m)),
\]
where $\bO_{\mu(m)}\subset \fg^\vee_D$ is the coadjoint orbit through $\mu(m)$. Due to the presence of
fixed points, however, this happens only on $M^o=(T^\vee\cP)^o$, and we have
\[
 \xymatrix{                                    & (T^\vee\cP)^o_{G,D,c}\ar[dr]^-{\mu}\ar[dl]_-{/G_D}& &\\
	      \Higgs^{rs}_{G,D,c}               &  & \fg_D^\vee\ar[r] & \fg^\vee_D \sslash G_D\\
}.
\]
Here $\fg_D^\vee\sslash G_D$ denotes, as usual, the GIT quotient, whose $\CC$-points correspond to closures of orbits.
We shall discuss  and refine this picture in the next section.

The Poisson structure, obtained by reduction from $(T^\vee\cP)^o$ turns out to coincide with the one defined in terms
of hypercohomology in (\ref{poisson_3}), see \cite{markman_thesis}, Corollary 7.15., and the integrability of the
former implies the integrability of the latter, which is defined on a bigger space.

      \section{Cameral covers and the Hitchin map}\label{cameral}
In this section we review  integrable system aspects of Higgs moduli. We start be recalling some
Lie-theoretic background, namely, the adjoint quotient morphism and its ``global analogue'' -- the Hitchin map. 
Then we turn to $L$-valued cameral covers and discuss very briefly generalised Prym varieties and abelianisation.
	       \subsection{Adjoint quotient and the Hitchin map}
		      \subsubsection{The Adjoint Quotient}
The group $G$ acts naturally on the coordinate ring of $\fg$ (via the coadjoint action) and by  a theorem of Chevalley (\cite{chevalley}) the subalgebra  of invariants
$\CC[\fg]^G\subset \CC[\fg]$ is isomorphic to a free algebra on  $l$ generators,
which are homogeneous of degree $d_j$, $1\leq j\leq l$.
 While the  choice of generators is largely non-unique, 
the set of degrees $\{d_j\}$ is  determined by $G$ and in fact  has a topological significance -- the Poincar\'e polynomial of
$G^{ad}$ is $p_{G^{ad}}(t)=\prod_{i}(1+ t^{2d_i-1})$.
The inclusion $\CC[\fg]^G\subset \CC[\fg]$ 
  corresponds to a  morphism (of affine varieties) $\chi: \fg\to \fg\sslash G=\spec \CC[\fg]^G$.
The $\CC$-points of the GIT quotient $\fg\sslash G$ are \emph{closures} of $G$-orbits in $\fg$.
 A specific   choice of $I_{j}\in \sym^{d_j}\fg^\vee$
with
$\CC[\fg]^G\simeq \CC[I_1,\ldots, I_l]$
determines an isomorphism $\CC^l\simeq \fg\sslash G$ and identifies $\chi$ with the  map
$x\mapsto (I_1(x),\ldots I_l(x))$.
Without the choice of $\{I_j\}$, $\fg\sslash G$ is not a vector space, but just  an affine
``cone'' -- a variety with a $\CC^\times$-action and a single fixed point.
The homothety action of $\CC^\times$ on $\fg$ descends to an action on $\fg\sslash G\simeq \CC^l$:
\[
 t\cdot (b_1,\ldots, b_l)= (t^{d_1}b_1,\ldots, t^{d_l}b_l)
\]
 and the morphism $\chi$ is  $\CC^\times$-equivariant. 
It is also $G$-invariant by construction.

Having chosen  Cartan and Borel subgroups $T\subset B\subset G$, we can use the embedding $W\subset GL(\ft)$ to gain another interpretation of $\chi$.  
Indeed, the inclusion $\CC[\ft]^W\subset \CC[\ft]$ gives rise to a quotient map $\ft\to \ft/W$, a finite flat morphism of affine varieties.
There exists a non-empty (Zariski-) open subset of $\ft$ -- the complement to the union of the root hyperplanes -- on which the quotient map is
an \'etale Galois cover with group $W$.
By a theorem   Chevalley,  the inclusion  $\ft\subset \fg$ determines an isomorphism  $\CC[\fg]^G\simeq \CC[\ft]^W$, and, consequently, 
$\ft/W\simeq \fg\sslash G$. We can think then of the adjoint quotient map as a morphism between affine varieties (in fact, affine spaces!)
$\fg\to \ft/W$.

Altogether, writing  $x^{ss}$ for  the semi-simple part of $x\in \fg$, we have the following descriptions of the morphism $\chi$:

\[
 \xymatrix{  \fg\ar[r]^-{\chi}   \ &\ \fg\sslash G         &\  \simeq       &  &\ft/W                              &\ \simeq          &\CC^l\\
	    x\ar@{|->}[r]        \ &\  \overline{G\cdot x} & \ar@{|->}[r]   &  &\left(G\cdot x^{ss}\right)\cap \ft & \ar@{|->}[r]     &(\underline{I}(x))\\
}.
\]
		  \subsubsection{The Hitchin map}
The adjoint quotient $\chi$ induces a morphism of  (total spaces of cone) bundles $\ad P= P\times_{Ad}\fg \to P\times (\fg\sslash G)$ which
can be twisted  with any $\CC^\times$-torsor, since $\chi$ is $\CC^\times$-equivariant. In particular, twisting with
$L^\times$ we get  a morphism of affine varieties
\[
\chi_{_{X,P}}: H^0(X, \ad P\otimes L)  \longrightarrow \cB= H^0(X, \ft\ctimes L/W)\simeq H^0(X, \bigoplus L^{d_i}).
\]
By the same token,
if 
 $T$ is a complex manifold and $\scP$  a holomorphic principal $G$-bundle on  $T\times X$,  
 $\chi$ induces a morphism 
\[
\chi_{_{T\times X,\scP}}:\ H^0(T\times X, \ad \scP\otimes p_2^\ast L) \longrightarrow H^0(T\times X, \left(\ft\ctimes p_2^\ast L\right)/W)= H^0(T,\cO_T)\otimes \cB.
 \]
Hence $\chi$ makes sense for families and
gives rise to a morphism, $h$,  from  $\Higgs_{G,D}$ to $\cB$ by 
 assigning to a $p_2^\ast L$-valued Higgs pair $(\scP,\Theta)$ on $T\times X$ the section
$\chi{_{T\times X,\scP}}(\Theta)\in H^0(T,\cO_T)\otimes \cB$, 
which is  nothing but a map  $T\to \cB$.
This morphism $h$ is the   \emph{Hitchin map}, given  on $\CC$-points by
$h([P,\theta])=\chi_{_{X,P}}(\theta)$, or, slightly informally, by
$h([P,\theta])= (I_1(\theta),\ldots, I_l(\theta))$, once the generators
$\{I_k\}$ are fixed.
 In what follows, we are going to suppress 
all  subscripts of $\chi$.

 The restrictions of $h$ to the respective connected components are known to be   \emph{proper} morphisms
\[
 h_c : \Higgs_{G,D,c}\to \cB
\]
which  endow $\Higgs_{G,D,c}$ with the structure of a Poisson ACIHS in the sense of \cite{donagi_markman}, Definition 2.9.
This means that $h_c$ is a proper flat morphism, which, away from  a closed subvariety $\Delta\subsetneq \cB$
has Lagrangian fibres, isomorphic to abelian varieties. The Lagrangian condition in the algebraic Poisson context is understood
generically: $Y\subset \Higgs_{G,D,c}$ is Lagrangian, if there is a symplectic leaf $S\subset \Higgs_{G,D,c}$, such that
$Y\subset \overline{S}$ and $Y\cap S\subset S$ is Lagrangian.  For the proof we refer to 
\cite{markman_thesis}, \cite{markman_sw} and \cite{bottacin}, as well as \cite{donagi_decomposition},
\cite{faltings}, \cite{scognamillo_elem} and
 \cite{don-gaits},
   extending the work in
\cite{hitchin_sd} and \cite{hitchin_sb}. The abelian varieties in question arise as generalised Prym varieties, associated to
cameral or spectral covers, see also (\ref{cameral_covers}) and (\ref{prym}).

The foliation of $\Higgs_{G,D,c}$ by closures of symplectic leaves is determined by the quotient $\cB/\cB_0$, where
\begin{equation}\label{B0}
 \cB_0 := H^0(X,\left(\ft\ctimes L/W\right)(-D))\simeq   H^0\left(X, \bigoplus_{i=1}^l L^{d_i}(-D)\right)\subset \cB.
\end{equation}
We have the following diagram, which  is a variant of \cite{markman_thesis}, Proposition 8.8:
\[
 \xymatrix{   & (T^\vee\cP)^{o}_{G,D,c}\ar[dr]^-{\mu}\ar[dl]_-{/G_D}& \\
	      \Higgs^{rs}_{G,D,c}\ar[d]_-{h_c}\ar[dr]^-{\overline{h}_c} & & \fg_D^\vee\ar[d]\\
		\cB\ar[r]& \cB/\cB_0\ar[r]^-{\simeq}&\fg^\vee_D \sslash G_D\\  
}
\]
and every $\overline{h}_c$-fibre contains a unique leaf of maximal rank.

	    \subsection{$L$-valued Cameral Covers}\label{cameral_covers}
The Hitchin base $\cB$   itself has modular interpretation -- it parametrises cameral covers of $X$, as we will review now.
The germ of the idea is already to be seen in Chevalley's theorem: a $G$-conjugacy class in $\fg$ can be identified with a $W$-conjugacy class
in $\ft$.
Twisting  the
(ramified)
 $W$-cover $\chi: \ft\to \ft/W$ with $L$ we get a
$W$-cover, $p$,  of the total space of $\ft\ctimes L/W\simeq \oplus_i L^{d_i}$
\[
 \xymatrix{  \tot \ \ft\ctimes L\ar[r]^-{p}\ar@{=}[d]& \tot \ft\ctimes L/W\ar@{=}[d]\\
	  \tot L^{\oplus l}\ar[r]&     \tot \bigoplus_i L^{d_i}\\
 }
\]
and pulling that back by the evaluation morphism 
 $\ev: \cB\times X\to \tot \ft\ctimes L /W$, $\ev(b,x)= b(x)$, we obtain a $W$-cover  $\cX=\ev^\ast p$ of $\cB\times X$. This is the 
\emph{universal cameral cover},  and each restriction of $\cX$ to $\{b\}\times X$
gives a cover $\pi_b: \widetilde{X}_b\to X$.

Our cameral curves $\widetilde{X}_b$ are all embedded  in  $\ft\ctimes L$ and inherit from it the 
$W$-action, which makes them
ramified Galois covers (with group $W$). It is worth mentioning that there is also an ``abstract'' notion of a cameral cover  -- one which 
does not involve the data of embedding into a vector bundle. Such covers are simply defined to be 
 locally (\'etale or analytically)
 the pullback of the cover $\chi$, see \cite{don-gaits}, Definition 2.6.

  The singularities and ramification behaviour of cameral covers are fairly well controlled. Indeed, $\chi$
 a singular hypersurface $\ft/W$,  the zero locus of the \emph{discriminant} 
\[
 \fD_\chi=\prod_{\alpha\in\cR}\alpha= (-1)^{|\cR|/2}\prod_{\alpha\in\cR^+}\alpha^2 =  P(\underline{I})\in \CC[\underline{I}]\simeq \CC[\ft]^W\subset \CC[\ft].
\]
The singular points of that hypersurface are $W$-orbits of semi-simple elements, lying on more than one  root hyperplane.
For instance, if $\fg=\fs\fl_3(\CC)$, the discriminant hypersurfaces is  a cuspidal cubic in $\ft/W\simeq \CC^2$, the cusp being the orbit of the
origin in $\ft\simeq \CC^2$.

Every root
$\alpha\in \ft^\vee$ gives  a morphism of bundles $\ft\ctimes L\to L$, which can be further pulled back to
$\tot \ft\ctimes L/W$ or $\tot \ft\ctimes L$. Consequently, the discriminant $\fD_\chi$ gives a morphism (as varieties over $X$) between
the total spaces of $\ft\ctimes L$ and $L^{|\cR|}$, and that can also be pulled further up
to $\fD\in H^0(\ft\ctimes L/W, q^\ast L^{|\cR|})$,
 as indicated on the next diagram

\[
\xymatrix{\widetilde{X}_b\ar@<-0.5ex>@{^{(}->}[r]\ar[d]^-{\pi_b}&\cX\ar[rr]^-{\iota}\ar[d]^-{\pi}& &\tot \ft\ctimes L\ar[d]_-{p}\ar[dr]^-{\widetilde{\fD}}&\\
	  \{b\}\times X \ar@<-0.5ex>@{^{(}->}[r]&\cB\times X\ar[rr]^-{\ev}\ar[dr]&& \tot \ft\ctimes L/W\ar[dl]_-{q}\ar[r]^-{\ev\fD}&\tot q^\ast L^{|\cR|}\\
	    & & X& &\\
}.
\]
 We denote by $Z(\fD)\subset \tot \ft\ctimes L/W$ the
vanishing locus of this section.

The possible singularities of $\widetilde{X}_b$ occur at the intersections of  root hyperplanes, i.e., over
points $b\in\cB$, where $b(X)$ meets the singular locus of  $Z(\fD)$.
We shall call a cameral cover \emph{generic} if it is smooth with simple Galois ramification, i.e., all ramification points have
ramification index one. That is, $b\in \cB$ is generic, if
 $\ev_b(X)\cap Z(\fD)^{sing}=\varnothing$ and $\ev_b(X)\pitchfork Z(\fD)^{sm}$.
We denote the locus of generic cameral covers by $\scB$, and
$\scB\varsubsetneq  \cB$ is a dense open subset.

	  \subsection{Generalised Pryms and abelianisation}\label{prym}
We recall here the definition of the generalised Prym variety, since the conventions used in the literature are
not completely uniform. As before, let $\Lambda_G:= \cchr_G\subset \ft$ be the cocharacter lattice. We have
that   $\Lambda_G \simeq \mhom(\CC^\times,T)$ and  $\Lambda_G \ztimes \CC^\times \simeq T$.
Donagi and Gaitsgory (\cite{don-gaits}) introduce two abelian sheaves, $\cT$ and $\overline{\cT}$, on $X$, associated
with the cover
 $\pi_o:\widetilde{X}_o\to X$.
The sections of the sheaf $\overline{\cT}$ on $U\subset X$ are $W$-equivariant (holomorphic) maps
$\pi_o^{-1}(U)\to T$, i.e., 
$\overline{\cT} =\pi_{o\ast} \left(\Lambda_G\otimes\cO^\times_{\widetilde{X}_o}\right)^W $. 
We note that the $W$-action on such maps incorporates both the $W$-action on $\widetilde{X}_o$ and
the $W$-action on $T$. In particular, an equivariant map must take the value $\pm 1$ on
root hyperplanes.
The sheaf
$\cT$ is the subsheaf of $\overline{\cT}$, whose sections over $U$ are the sections in $\overline{\cT}(U)$, taking value
$+1$ on $\widetilde{X}_o\cap \{\alpha=0\}$ for all roots $\alpha\in\cR$.
Since we are assuming that $G$ is simple, 
 one can see (\cite{don-pan}, Lemma 3.3) that
  $\cT= \overline{\cT}$ if $G\neq B_l$. In the exceptional case, 
$\overline{\cT}/\cT$ is   $\ZZ/2\ZZ$-torsion, supported at the branch points of $\pi_o$.
Moreover, \emph{ibid.}, Claim 3.5,   $H^1(X,\cT)$ and $H^1(X,\overline{\cT})$ are isogenous abelian varieties.
The 
\emph{generalised Prym variety } associated to the given cameral cover is  $\Prym_{\widetilde{X}_o/X}:=H^1(X,\cT)$.
These varieties are isogenous to 
$\left( \Lambda_G\ztimes H^1(\widetilde{X}_o,\cO^\times) \right)^W$, which is the set of $W$-invariant $T$-bundles
on $\widetilde{X}_o$.

We denote by $\Prym_{\cX/\scB}$ the relative Prym fibration (over $\scB$) associated with $\cX$.
By the abelianisation theorem (\cite{don-gaits}), $h_c^{-1}(o)$ is a  $\Prym_{\widetilde{X}_o/X}^0$-torsor, and, moreover,
$\Higgs_{G,D,c}$ is a $\Prym^0_{\cX/\scB}$-torsor. The two can be locally identified by choosing local sections (over $\scB$).

We remark that Donagi and Gaitsgory describe explicitly  
spectral data, corresponding to
$h_c^{-1}(o)$, i.e., the particular $\Prym^0_{\widetilde{X}_o/X}$-torsor, see \cite{don-gaits}, Theorem 6.4 and 
\cite{don-pan}, Appendix A.1. We do not need this description in what follows, so will not dwell on it.
      \section{The Infinitesimal Period Map}\label{cubic}
	  \subsection{The main theorem}
We now have at hand all ingredients needed for stating  the  main result.
     Fix a point   $o\in \scB$ and a topological type $c\in \pi_1(G)$, such that 
$\Higgs_{G,D,c}\neq \varnothing$.
The base point  corresponds to a maximal rank symplectic leaf $S$, whose closure in $\Higgs_{G,D,c}$
is $h_c^{-1}(\{o\} + \cB_0)$. In general, this closure is strictly bigger than the one in Markman's construction, which takes
place on the smaller locus $\Higgs_{G,D,c}^o$ of Higgs pairs,  having a regularly stable underlying bundle. The difference, however, 
is away from the generic locus of $\cB$, to which we now restrict.
We consider the set
\begin{equation}\label{boldB}
 \bB:= \left(\{o\}+ \cB_0\right)\cap \scB  \subset \cB
\end{equation}
which supports an integrable system  (in the symplectic sense), all of whose fibres are proper:
\[
 \xymatrix@1{  \Higgs_{G,D,c} & \supset  &\left. \cS\right|_\bB= h_c^{-1}\left(\bB \right)\ \ar[r]^--{h_\bB} &\ \bB \ni o&}
\]
where $h_{\bB}=\left. h_{c}\right|_{\bB}$. Our main theorem is a statement about the infinitesimal period map of this family
of abelian torsors.

\begin{thma}[\cite{ugo_peter_cubic}]
There exists a natural isomorphism 
\[
 T_{\bB,o}\simeq H^0(\widetilde{X}_o,\ft\ctimes K_{\widetilde{X}_o})^W .
\]
Let $Y_\xi \in T_{\bB,o}$ denote both the preimage of  $\xi \in H^0(\ft\ctimes K_{\widetilde{X}_o})^W $  and the
corresponding  constant vector field.
Under this isomorphism,    
the differential at $o\in \bB$ of the  period map of  $h_{\bB}: \left. \cS\right|_{\bB}\to \bB$
is given by  
\[
 c_o: H^0(\widetilde{X}_o, \ft\ctimes K_{\widetilde{X}_o})^W \longrightarrow \sym^2 \left(H^0(\widetilde{X}_o, \ft\ctimes K_{\widetilde{X}_o})^W \right)^\vee
\]
\[
 c_o(\xi)(\eta,\zeta)= \frac{1}{2}\sum_{p\in\textrm{Ram}(\pi_o)} \textrm{Res}_p^2\ \left( \pi_o^\ast \left.\frac{\cL_{Y_\xi}(\fD)}{\fD}\right|_{\{o\}\times X}\eta\cup\zeta\right).
\]

\end{thma}
    	  \subsection{Related results}
We should note that several instances of {\bf Theorem A} (and its reformulation,  {\bf Theorem B} in \cite{ugo_peter_cubic}) have already been 
established  in the literature.
First of all, for the usual Hitchin system ($D=0$) and $G=SL_n(\CC)$ the formula appears in unpublished work of T.Pantev,
while for the case $G=SL_2(\CC)$ the formula can be found in \cite{geom_trans}, (47). Building on that, D.Balduzzi
(\cite{balduzzi}) dealt with the case of a semi-simple structure group (still in the case  $D=0$). The same formula for the cubic is obtained
(by a somewhat different method) in \cite{hertling_hoev_posthum}. When it comes to the  generalised Hitchin system,
the only similar result that we are aware of appears in the context 
 of the Neumann oscillator, see  \cite{hoevenaars}. There 
the base curve is $\PP^1$, while $D=2\cdot \infty +\sum_{i} n_i q_i$  and $G=SL_2(\CC)$. 
    \subsection{Proof of the theorem}
We break the argument into several steps.
	\subsubsection{Step 1: the isomorphism} 
Since $\bB$ is an open subset of an affine space modelled on $\cB_0$,  we have
a canoncal isomorphism $T_{\bB,o}= \cB_0$. This isomorphism, however, makes no reference to $o\in \bB$.
So we  restrict  $\cX$ to $\bB\times X$, and consider the embedding 
 $\widetilde{X}_o\hookr\left. \cX\right|_{\bB\times X}$. Setting
$N$ to be  the normal bundle of the inclusion 
$\io_o \widetilde{X}_o\subset \tot\ft\ctimes L$ and  $r:\tot \ft\ctimes L\to X$ to be the bundle projection, 
we get an isomorphism
$ T_{\bB,o} \simeq   H^0(\widetilde{X}_o, N(-r^\ast D))^W $. 
 We now claim that there is a natural isomorphism
\[
 H^0(\widetilde{X}_o, N(-r^\ast D))^W\simeq H^0(\widetilde{X}_o, \ft\ctimes K_{\widetilde{X}_o})^W
\]
For that, recall that for any vector space $V$ there is a natural $\ft$-valued 2-form on $V\oplus \left(V^\vee\otimes \ft\right)$, namely
\[
\left((x,\alpha\otimes s),(y,\beta\otimes t)\right)= \alpha(y)s-\beta(x)t.
\]
Here we could have replaced the Cartan subalgebra $\ft$ by any vector space -- for instance, replacing it by $\CC$ would give the
canonical symplectic form on $V\oplus V^\vee$. This form induces a $\ft$-valued 2-form on $\tot \ft\ctimes K_X$, and, after twisting with
$D$, a meromorphic 
$\ft$-valued 2-form
$\omega_\ft\in H^0(\tot \ft\ctimes L, \ft\ctimes \Omega_{\ft\ctimes L}^2(r^\ast D))^W$.
Contraction with $\omega_\ft$  gives a 
sheaf homomorphism
\[
  N\to \ft\ctimes K_{\widetilde{X}_o}(r^\ast D).
\]
 While it is not necessarily an isomorphism,  it does induce an isomorphism on  invariant global sections. This follows from
\cite{kjiri}, where  (following the reasoning in  \cite{hurtubise_kk}), the author shows that the generalised Hitchin system
satisfies the rank-2 condition of Hurtubise and Markman. That turns the statement into a special case of Proposition 2.11 in 
\cite{hurtubise_markman_rk2}.
	  \subsubsection{Step 2: Recasting the symplectic structure }
The relative Prym fibration
$g_{\bB}: \Prym_{\cX/\bB}^0\longrightarrow \bB$ is Lagrangian
and  under the local identifications
$\left. \cS\right|_{\cU}\simeq \Prym^0_{\cX/\cU}$, $\cU\subset \bB$,  the  symplectic structures on both sides coincide.
Indeed, for the abelian variety   $P_o=\Prym_{\widetilde{X}_o/X}^0$ we have
\[
 T_{P_o}=H^1(\widetilde{X}_o,\ft\ctimes\cO)^{W}\ctimes \cO_{P_o},
\]
and so,  by Serre duality,  for any $\scL\in P_o$ we get
\[
 T_{P_o, \scL}^\vee = T_{\bB,o}\simeq H^0(\widetilde{X}_o, \ft\ctimes K_{\widetilde{X}_o})^W.
\]
This gives the Lagrangian structure on the Prym fibration (restricted to $\bB$), and it coincides with the one obtained by Marsden--Weinstein reduction.
For a very concrete description in the case $X=\PP^1$, see Theorem 1.10 in \cite{hurtubise_kk}.
	\subsubsection{Step 3: Reduction to a Kodaira--Spencer calculation}

Any  family $h: \cH\to \bB$ of polarised compact (connected) K\"ahler manifolds gives rise to  a weight-1 
polarised $\ZZ$-VHS 
 $\left(\cF^\bullet,   \cF_\ZZ, \nabla^{GM}, S\right)$
on $\bB$ with  a period map $\Phi$. By a theorem of Griffiths (\cite{griffiths_periods}, part II, Theorem 1.27) we have that
$d\Phi_o= m^\vee\circ \kappa$, where $\kappa: T_{\bB,o}\to H^1(\cH_o, T)$ is the Kodaira--Spencer map and
\[
 m^\vee: H^1(\cH_o, T)\to H^1(\cH_o, \cO)\otimes H^0(\cH_o, \Omega^1)^\vee\simeq_S \left(\cF^{1\vee}_o\right)^{\otimes 2}
\]
is induced  by cup product $H^1(T)\times  H^0(\Omega^1)\to H^1(\cO)$.  

By a choice of local section we can replace the family $\left. \cS\right|_\bB$ by $\cH=\Prym_{\cX_\bB/\bB}$ and 
from \emph{Step 2} we have  that
$H^1(P_o, T_{P_o})\simeq H^1(\widetilde{X}_o,\ft\ctimes\cO)^{W\otimes 2}$.
 Polarisation-preserving deformations
are contained in
 \[
\sym^2 H^1(\widetilde{X}_o,\ft\ctimes\cO)^{W}\simeq \sym^2 H^0(\widetilde{X}_o, \ft\ctimes K_{\widetilde{X}_o})^{W \vee}.
\]
But by properties of cup product,  $m^\vee$ is dual to the 
multiplication map
\[
m: H^0(\widetilde{X}_o, \ft\ctimes K_{\widetilde{X}_o})^{W \otimes 2}\to H^0(\widetilde{X}_o,\ft\otimes \ft \otimes K^2)\to_{\textrm{tr}} H^0(\widetilde{X}_0,K^2).
\]
In this way, all data are expressed in  terms of the family of cameral curves
$f=\textrm{p}_1\circ \pi: \left.\cX\right|_{\bB}\to \bB$, and the Kodaira--Spencer maps
coincide, since
$g_{\bB \ast}T_{{\Prym/\bB}}\simeq f_\ast T_{\cX/\bB}$. 
Moreover, 
 the polarisation on the Pryms is determined by the polarisation on the $\widetilde{X}_b$. 
Hence
the question of computing the infinitesimal period map of
$\cS_\bB\to \bB$ is replaced with  the same question, but for the family 
$\cX_\bB\to \bB$.

Finally, since for a finite dimensional vector space $V$  the  natural isomorphism
 $\mhom(V^\vee,\mhom(V^\vee,V))=\mhom(V^{\vee\otimes 3},\CC)$ is given by
\[
F \longmapsto \left(Y\otimes \alpha\otimes\beta\mapsto \beta(F(Y)(\alpha))\right), 
\]
we obtain the following 
\begin{proposition} 
The differential of the period map of $h_\bB: \left. \cS\right|_{\bB}\to \bB$ at $o$
is given by
\[
 c_o: H^0(\widetilde{X}_o, \ft\ctimes K_{\widetilde{X}_o})^W \longrightarrow \sym^2 \left(H^0(\widetilde{X}_o, \ft\ctimes K_{\widetilde{X}_o})^W \right)^\vee,
\]
\[
 c_o(\xi)(\eta,\zeta)=\frac{1}{2\pi i}\int_{\widetilde{X}_o}\kappa(Y_\xi)\cup \eta\cup \zeta.
\]
where $\kappa$
is the Kodaira--Spencer map of the family $\left. \cX\right|_{\bB}\to \bB$ at $o\in\bB$.
\end{proposition}
	\subsubsection{Step 4: Kodaira--Spencer calculation}
If one is given a family of compact K\"ahler manifolds over a contractible base, then
any holomorphic vector field on the base can be lifted to a smooth vector field on the total space of
the family. In general, however, there does not exist a holomorphic lift. Locally, such lifts do
exist, and the Kodaira--Spencer map measures the obstruction to the existence of a global one.
More intrinsically, one can describe this map as a connecting homomorphism for the
derived image of the projection morphism.

Indeed, getting back to our setup, we see that pushing forward the tangent sequence  of $f: \left.\cX\right|_{\bB}\to \bB$ gives  a connecting homomorphism
$\delta: T_\bB \to R^1 f_\ast T_{\cX_\bB/\bB}$.
For a contractible neighbourhood $\cU\subset \bB$ of $o\in \bB$, the (global) Kodaira--Spencer map $\varkappa$ over $\cU$ is 
the  induced map  on global sections:
\[
\varkappa=\Gamma_\cU(\delta):\ \Gamma_\cU (T_\bB)\to \Gamma_\cU( R^1 f_\ast T_{\cX_\bB/\bB})= H^1(\left. \cX\right|_{\cU}, T_{\cX_\cU/\cU}),
\]
while the (pointwise) Kodaira--Spencer map   appearing in Griffiths' theorem is 
$\kappa = \varkappa_o: T_{\cU,o}\to H^1(\widetilde{X}_o,T_{\widetilde{X}_o})$, 
obtained by passing to the fibre over
$o$.

We  compute $\varkappa$  on  a convenient \v{C}ech cover of
$\cX_\bB$, thus following  the original approach of  Kodaira and Spencer, see \cite{kodaira-spencer},
or \cite{kodaira_defo}, Ch.4. Since we have a family of covers of a fixed target curve, $X$, we choose a covering which 
facilitates the
 handling of  varying  branch loci, as follows.

Let $\ram(\pi)$ and $\bra(\pi)$ be the ramification and branch loci of $\pi: \cX_\cU\to \cU\times X$.
We set $\bU:= \cX_\cU\backslash \ram(\pi)$   and introduce the 
  cover $\cX_\cU=\bU\cup\bV$, where $\bV\supset \ram(\pi)$ is a certain tubular neighbourhood,
 constructed as 
follows.  Consider $\bra(\pi_o)=\{p_1,\ldots,p_N\}\subset X$, $N=|\cR|\deg L$,  and  choose an atlas
$\left\{(U_j,z_j), j=0\ldots N  \right\}$ of $X$,  where 
$U_0=X\backslash\bra(\pi_o)$, and $\{U_j\ni p_j\}$  non-intersecting open disks. For simplicity, we assume
$\textrm{supp}(D)\cap \bra(\pi_o)=\varnothing$.
Since $\cU\subset \scB$, by the genericity assumption
$ \cU\times X\supset \bra(\pi)_\cU\to \cU$ is an unramified $N:1$ cover, and admits, by the implicit function theorem,  local sections $c_j:\cU\to X$,
$1\leq j\leq N$, 
such that $c_j(o)=p_j$, and  $c_j(\cU)\subset U_j$ (possibly after shrinking $\cU$). We then  define
\[
\bV:= \pi^{-1}\left( \coprod_{j\neq 0}\textrm{graph }c_j \right)\subset \cX_\cU.
\]
This set has $\deg L |\cR| |W| /2$ connected components, which we  index as $\bV_{j\alpha}^k$, and group them into
$\bV_{j\alpha}=\coprod_k \bV_{j\alpha}^k=\pi^{-1}\left(\textrm{graph }c_j\right)\cap \{Z(\alpha)\}$,
$\bV=\coprod \bV_{j\alpha}$. Let us note that the cover $\{\bU,\bV\}$ is good and that for the calculation of \v{C}ech
1-cochains we have to focus attention on $\pi^{-1}\left(\cU\times \coprod_{j\neq 0} U_j\right)\supset \bU\cap \bV$.
We have local equations for  $\cX_{\cU\times U_j}$ of the form
\[
 \left|
	  \begin{array}{l}
	   I_1(\alpha_1,\ldots,\alpha_l)=b_1(\underline{\beta},z)\\
	    \ldots \\
	   I_l(\alpha_1,\ldots,\alpha_l)=b_l(\underline{\beta},z)\\
	  \end{array}
\right. ,
\]
where  $\underline{\beta}$ are   coordinates on $\cU$, say obtained from a choice of basis on $\cB_0$, 
and $z=z_j$.
On every   component of $\bU$, we can use as (\'etale!) coordinates $(\underline{\beta},z)$, giving a parametrisation
$\alpha_i= g_i^0(\underline{\beta}, z)$ (after passing to a cover). This can then be carried to the other components by the $W$-action, 
$s_j\cdot g_i = g_i- n_{ij}g_j$.
On the other hand, since $\textrm{supp} D\subset \left(X\backslash \coprod_{j\neq 0}U_j\right)$,   the local trivialisations
of $K_X$ (provided by the atlas)  turn roots into maps $\ft\ctimes L_{U_j}\to \CC$.
By genericity and WPT on $\bV^k_{j\alpha}$, say $\alpha:=\alpha_1$, we have
\[
 \left|
      \begin{array}{lr}
       \alpha^2 = (z-c(\underline{\beta}))v(\underline{\beta},z)& \\
       \alpha_i=g_i(\alpha,z),& i\geq 2\\
      \end{array}
\right. .
\]
Here the holomorphic function $v$ is non-zero along $\textrm{graph}(c)\subset \cU\times U_j$. Inverting it, we have
 $z=\phi(\underline{\beta},\alpha) =c(\underline{\beta}) + \alpha^2 u(\underline{\beta},\alpha)$, for some other holomorphic function $u$.

Starting with the (constant)  vector field $Y= \partial_\beta\in \Gamma_{\cU}(T_\bB)$, the
 cocycle $\{\varkappa_{\alpha z}(Y)\}$ representing $\varkappa(Y)$ is given, on $\bV_{j\alpha}^k\cap\bU$ by the
vertical vector field
\[
\left. \frac{\partial \alpha}{\partial \beta}\right|_{z=\phi(\beta,\alpha)}\frac{\partial}{\partial \alpha}.
\]
By the implicit function theorem, 
\[
 \left. \frac{\partial\alpha}{\partial \beta}\right|_{z=\phi(\underline{\beta},\alpha)}= \frac{\alpha}{2}
\left. \frac{\partial_\beta \alpha^2}{\alpha^2}\right|_{z=\phi(\underline{\beta},\alpha)}= 
-\frac{\partial_\beta c}{2\alpha u(\underline{\beta},\alpha)} + \left. \frac{\alpha}{2}\partial_\beta v\right|_{z=\phi(\underline{\beta},\alpha)}.
\]

Notice that  along $\alpha=0$ the first term has a pole, while the second has a zero. On the other hand, the
``logarithmic derivative'' of $\fD$ along $Y$ is
\[
 \pi^\ast \frac{\partial_\beta \fD}{\fD}=\sum_{\alpha\in\cR^+} \frac{\partial_\beta \alpha^2}{\alpha^2} 
\]
and so, up to a factor of $1/2$, taking quadratic residues with respect to the latter has the same effect as contracting
 with $\varkappa_{\alpha z}(Y)$ and taking residues. Summing over all branch points, roots, and $\ZZ/2\ZZ$ cosets in $W$ completes the proof.

\qed

      \section{Some Related Geometries}\label{survey}
	  \subsection{The Donagi--Markman cubic condition}
In this section we recall very briefly some structures, naturally related to the
study of the infinitesimal period map for algebraic integrable systems. We also  give references to 
   recent developments in these areas.

As discussed in the introduction, many of the  differences between real and algebraic (analytic) integrable systems
are rooted in the fact that complex tori have moduli. One such striking difference arises when one considers  
 a holomorphic family, $h: \cH\to \bB$, of complex tori,
with the property that $2\dim \bB= \dim \cH$. In the smooth category, locally on $\bB$ there is always a 
compatible Lagrangian structure on the family, i.e., symplectic structure, for which the fibres are Lagrangian.
In the holomorphic or algebraic context, however, there is a local obstruction to the existence of such structure, as discovered by
Donagi and Markman. The presence of Lagrangian structure forces the infinitesimal period map to be a section of
$\sym^3 T^\vee_\bB$, rather than just  a section of $T_\bB^\vee\otimes \sym^2 T^\vee_\bB$. After making the appropriate choices, this condition
forces the period matrix to be locally a Hessian of a holomorphic function, known as \emph{holomorphic prepotential} $\cF:\cU\to \CC$.

We sketch now how this obstruction arises. Consider a contractible open set $\cU\subset \CC^d$ and a holomorphic map 
$\Phi: \cU\to \HH^{d}\subset \textrm{End}(\CC^{d})$, 
where $\HH^{d}$ is Siegel's upper-half space of dimension $d$. Let also $\Omega =(\left. 1 \ \right| \Phi): \cU\to \textrm{Hom}(\CC^{2d},\CC^{d})$
be the map $s\mapsto \Omega(s)= (\left. 1 \right| \Phi(s))$. Finally,  let $\Gamma\simeq \ZZ^{2d}$ be the group
 of holomorphic automorphisms  of $\cU\times \CC^d$, generated by
$(s,z)\mapsto (s, z+ \Omega(s)(e_j))$, $j=1\ldots 2d$. 
Then the quotient  $\cH= \cU\times\CC^d/\Gamma$ is a family (over $\cU$) of abelian varieties.
We  can happily endow $\cU\times \CC^d\simeq T^\vee_\cU$ with the canonical symplectic structure,
but this structure will  not descend to the quotient, unless the sections $s\mapsto \left(s, \Omega(s)(e_j)\right)$ of $T^\vee_\cU$,
determined by the columns $\Omega(e_j)$, happen to be \emph{Lagrangian}.
This happens precisely when the 1-forms $\sum_k \Omega_{kj}ds_k$ are \emph{closed}, i.e., 
$\partial_i\Omega_{kj}=\partial_k\Omega_{ij}$. The latter
 implies that
  $\Phi=\textrm{Hess }\cF$ for some holomorphic function $\cF$, possibly after shrinking $\cU$.
Correspondingly,  $c= d\Phi = \sum \partial^3_{ijk}\cF dx_i\cdot dx_j\cdot dx_k$.

We direct the reader to the beautiful expositions in
\cite{donagi_markman_cubic}, \S 1 and \cite{donagi_markman}, \S 7 for a different version of this argument,  discussion and applications.
For the case of the generalised Hitchin system, our
 {\bf Theorem B} in \cite{ugo_peter_cubic} contains a formula for the infinitesimal period map
which makes it evident that $c$ is a
section of $\sym^3T^\vee_\bB$. 
	  \subsection{Special K\"ahler  and Seiberg--Witten geometry}

The data of an  ACIHS $h:\cH\to \bB$ naturally gives rise to a certain kind of K\"ahler geometry on $\bB$. Indeed, 
the imaginary part of the period matrix $\textrm{Hess }\cF$ is symmetric and positive-definite, and hence can be used to
define a K\"ahler metric. 
 More intrinsically, 
away from the discriminant locus we have an   identification between
the vertical bundle (the direct image of the relative tangent bundle of $h$) and the cotangent bundle to $\bB$. 
Choosing a local section of $h$ over $\cU\subset \bB$ we can identify $\cH\to \bB$ with its relative Albanese fibration,
i.e., identify $\cH_\cU\to \cU$ with a family of polarised abelian varieties.
 The polarisation gives rise to a translation-invariant metric on each fibre, and we obtain a ``semi-flat'' metric on $\bB$.
Furthermore,  the  bundle of lattices can be used to define a flat connection on
$T^\vee_\bB$, and hence on $T_\bB$.

 This kind of  differential-geometric structure is known as  \emph{special K\"ahler geometry}. Abstractly,
one starts with a K\"ahler manifold $(M,I,\omega)$, $\bB=(M,I)$, and considers 
  a flat, symplectic, torsion-free
connection $\nabla$ (on $T_{M}$), such that $d^\nabla I=0$.
 The special K\"ahler geometry, arising from an integrable system carries an additional integrality property.
Such structures were introduced by the physicists (\cite{bcov}) for the study of vector multiplets in four-dimensional
$N=2$ SUSY. For a beautiful -- and by now, classical -- introduction to these structures, one can consult
\cite{danfreed}. 
The generalised Hitchin system depends on a large amount of input data, including  the genus of $X$, the degree of the divisor
$D$ and the group $G$. If one wants build realistic  physical examples these discrete invariants are quite constrained.
 Donagi and Witten 
have discussed various proposals of this kind and made some tests in \cite{donagi_witten}. 
For a discussion of the relations between Seiberg--Witten theory, integrable systems in general and
the generalised Hitchin system in particular, one can consult
\cite{donagi_swis} and \cite{markman_sw}. For some recent developments one can consult  the other articles in this
volume.

	  \subsection{VHS and $tt^\ast$-geometry}
In the proof of {\bf Theorem A} we have used explicitly several bits of elementary Hodge theory, and we have, on the other hand, 
indicated that the presence of an ACIHS on $\bB$ is equivalent to the presence of an (integral) special K\"ahler structure.
This is by no means  a coincidence:  one can show   (\cite{hertling_ttstar}, \cite{bartocci_mencattini_sk}) that the
special K\"ahler data $(M,I,\omega,\nabla)$  is equivalent to the data of a weight-1 real polarised variation of Hodge structures on
$T_{\bB,\CC}$. The latter  means that we have a polarised, weight one, real VHS on $\bB$,
with Hodge flag  $\cF^1= T_\bB\subset \cF^0$, where $\cF^0\otimes_{\cO_\bB}\scC^\infty_\bB= T_{\bB,\CC}$. This forces
$\cF^0$ to be a
(holomorphic) extension of $T^\vee_\bB$ by $T_\bB$. To construct this extension from  the special K\"ahler connection,
one considers the type decomposition
$\nabla=\nabla'+\nabla''$ and sets $\dbar_{\cF^0}=\nabla''$. The  
 Gauss--Manin connection of the VHS is the  set  to be 
$\nabla^{GM}=\nabla'$, the $(1,0)$-part of $\nabla$.

 The  kind of VHS appearing above falls within the context
considered initially by C.Simpson in \cite{simpson_uniformisation}.
Indeed,  the
second fundamental form of the Gauss--Manin connection $\nabla^{GM}: \cF^0\to \cF^0\otimes T^\vee_\bB$
is the induced $\cO_B$-linear map $\cF^1\to (\cF^0/\cF^1)\otimes T^\vee_\bB$.
This can be identified with 
a $T_\bB^\vee$-valued Higgs field $\theta$ on the
associated graded bundle $\textrm{gr }\cF^\bullet = T_\bB\oplus T_{\bB}^\vee$. This is a nilpotent Higgs field, whose  only (possibly)  non-zero component is 
contained in
$H^0(\bB, T_\bB^{\vee\otimes 3})$. This  is precisely the Donagi--Markman cubic of the corresponding integrable system. Considering
the Hermite--Yang--Mills equation for the Higgs vector bundle $(T_\bB\oplus T_\bB^\vee, \theta)$ 
one finds  the $tt^\ast$-equation (\cite{danfreed}, 1.32)
which is  one of the main objects of interest in the original work of 
\cite{bcov}.

Thus, starting from an ACIHS, and considering  its associated special K\"ahler geometry as a particular
variation of Hodge structures on the base $\bB$, one stumbles upon a piece of non-abelian Hodge theory:
a Higgs bundle on $T_\bB\oplus T_\bB^\vee$, for which the Hitchin--Kobayashi correspondence is tied
with $tt^\ast$-geometry. These phenomena seem to be rooted in  mirror symmetry and non-commutative geometry.

	  \subsection{Bryant--Griffiths geometry and non-commutative Hodge structures}
The relation between Hodge theory, $tt^\ast$-geometry and mirror symmetry is
 extremely intricate and far from being  completely understood. Here we  indicate  only the main points and give references
to the some of the literature.

 One of the formulations of the mirror conjecture is in terms of Frobenius manifolds. On the
$A$-side the Frobenius structure is obtained from the quantum cohomology product, while on the $B$-side it is phrased in terms
of various extensions of the notion of Hodge structure. These  extended structures seem to be  most conveniently encapsulated in the
  notion of \emph{non-commutative} Hodge structure, see \cite{kkp_hodge_ms}, \cite{sabbah_nch}. Variations of such
structures on tangent bundles of manifolds give rise to $tt^\ast$-geometry. 

When  applied to the context of  mirror symmetry of Calabi--Yau threefolds, the $B$-side Frobenius geometry can be
expressed in elementary terms and related to  Bryant--Griffiths (\cite{bryant_griffiths}) type geometry.
We recall that the Bryant--Griffiths setup is the defining example of \emph{projective} special K\"ahler geometry
($N=2$ supergravity in physics), where the weight-one $\RR$-VHS on $T_{\bB,\CC}$ is refined to a weight-three VHS
(\cite{danfreed}, \S 4).
 This relation between Frobenius-type structures and projective special K\"ahler geometry is discussed in great detail
in \cite{hertling_hoev_posthum}, but with a somewhat idiosyncratic nomenclature, originating in \cite{hertling_ttstar}.

	  \subsection{Large $N$ duality and ADE Hitchin systems}
It turns out that if $G$ is a group of  $ADE$ type,  then the  base $\cB_\fg$ of the usual ($D=0$) Hitchin system does support the 
kind of Frobenius-like structures mentioned above.
The relation between Hitchin systems and moduli of Calabi--Yau threefolds was discovered in
 \cite{geom_trans} and \cite{ddp}.
 There the authors  construct a family $\cX\to \cB_\fg$ of surface-fibred quasi-projective Calabi--Yau threefolds, 
for which the family of  intermediate Jacobians is isogenous to $\Prym^0_{\cX/\cB}$ and the Yukawa cubic is identified with the
Donagi--Markman cubic.

The family of 3-folds $\cX\to \cB_\fg$ arises as a family of surfaces over $\cB_\fg\times X$, and  is constructed as follows.
Let $\fg$ be a simple Lie algebra of type ADE, and
 $\Gamma\subset SL_2(\CC)$  a finite subgroup, corresponding to $\fg$ by the McKay correspondence. Let $X$ be a smooth curve of genus $g\geq 2$
and $V\to X$ a $\Gamma$-equivariant rank two vector bundle, with a fixed isomorphism $\det V\simeq K_X$.
To this data one associates the Calabi--Yau threefold $Y_0= \tot (V)/\Gamma\to X$, all of whose fibres (over $X$) are ALE spaces of type
$\Gamma$. 
Then, by work of B. Szendr\"oi (\cite{szendroi_artin}, \cite{szendroi_3folds}), there exists a family of
surfaces $\cQ\to \tot \ft\ctimes K_X/W$, uniquely characterised by two properties: that its restriction to
 the zero section (of $\ft\ctimes K_X/W$) is isomorphic to $Y_0$ and that its restriction to a fibre
  is isomorphic to the universal unfolding of $\CC^2/\Gamma$. Then the family of threefolds is defined as
$\cX:= \ev^\ast \cQ$, where, as
in \S \ref{cameral_covers}, $\ev:\cB_\fg\times X\to \tot \ft\ctimes K_X/W$ is the evaluation morphism.

  In \cite{hertling_hoev_posthum} the corresponding Frobenius-type structures on $\cB_\fg$ are discussed.

We  hope that our understanding of the Donagi--Markman cubic  may  be useful in unravelling some of the
intricacies of the analogous story when $D\neq 0$. Judging by the recent development in \cite{cddp}, however, that will require
developing a considerable  amount of new techniques.

\bibliographystyle{alpha}
\bibliography{biblio.bib}

\newcommand{\etalchar}[1]{$^{#1}$}
\begin{thebibliography}{DDD{\etalchar{+}}06}

\bibitem[AAB00]{biswas_anchouche_azad_prb}
B.~Anchouche, H.~Azad, and I.~Biswas.
\newblock Holomorphic principal bundles over a compact {K}\"ahler manifold.
\newblock {\em C. R. Acad. Sci. Paris S\'er. I Math.}, 330(2):109--114, 2000.

\bibitem[AB01]{biswas_anchouche_eh}
B.~Anchouche and I.~Biswas.
\newblock Einstein-{H}ermitian connections on polystable principal bundles over
  a compact {K}\"ahler manifold.
\newblock {\em Amer. J. Math.}, 123(2):207--228, 2001.

\bibitem[AG88]{arnold_givental}
V.I. Arnold and A.B. Givental'.
\newblock Symplectic geometry.
\newblock In {\em Dynamical systems {I}{V}}, volume~4 of {\em EMS}, pages
  1--136. Berlin Heidelber Springer: New York, 1988.

\bibitem[Ati57]{atiyah_elliptic}
M.~F. Atiyah.
\newblock Vector bundles over an elliptic curve.
\newblock {\em Proc. London Math. Soc.}, 3(7):414--452, 1957.

\bibitem[Bal06]{balduzzi}
D.~Balduzzi.
\newblock Donagi-{M}arkman cubic for {H}itchin systems.
\newblock {\em Math. Res. Lett.}, 13(5-6):923--933, 2006.

\bibitem[BBNN06]{bal_bis_nag_news}
V.~Balaji, I.~Biswas, D.~S. Nagaraj, and P.~E. Newstead.
\newblock Universal families on moduli spaces of principal bundles on curves.
\newblock {\em Int. Math. Res. Not.}, pages Art. ID 80641, 16, 2006.

\bibitem[BCOV94]{bcov}
M.~Bershadsky, S.~Cecotti, H.~Ooguri, and C.~Vafa.
\newblock Kodaira-{S}pencer theory of gravity and exact results for quantum
  string amplitudes.
\newblock {\em Commun. {M}ath.{P}hys.}, (165):311--428, 1994.

\bibitem[BD14]{ugo_peter_cubic}
U.~Bruzzo and P.~Dalakov.
\newblock Donagi-{M}arkman cubic for the generalized {H}itchin system.
\newblock {\em Internat. J. Math.}, 25(2):1450016, 20, 2014.
\newblock \href{http://xxx.lanl.gov/abs/1308.6788}{arXiv 1308.6788}.

\bibitem[BG83]{bryant_griffiths}
R~Bryant and Ph~Griffiths.
\newblock Some observations on the infinitesimal period relations for regular
  threefolds with trivial canonical bundle.
\newblock Arithmetic and geometry, Pap. dedic. I. R. Shafarevich, Vol. II:
  Geometry, Prog. Math. 36, 77-102 (1983)., 1983.

\bibitem[BGO11]{ugo_bea_ss}
U.~Bruzzo and B.~Gra{\~n}a~Otero.
\newblock Semistable and numerically effective principal ({H}iggs) bundles.
\newblock {\em Adv. Math.}, 226(4):3655--3676, 2011.

\bibitem[BH12a]{biswas_hoffmann_poincare}
I.~Biswas and N.~Hoffmann.
\newblock Poincar\'e families of {$G$}-bundles on a curve.
\newblock {\em Math. Ann.}, 352(1):133--154, 2012.

\bibitem[BH12b]{biswas_hoffmann_torelli_bun}
I.~Biswas and N.~Hoffmann.
\newblock A {T}orelli theorem for moduli spaces of principal bundles over a
  curve.
\newblock {\em Ann. Inst. Fourier (Grenoble)}, 62(1):87--106, 2012.

\bibitem[BM09]{bartocci_mencattini_sk}
C.~Bartocci and I.~Mencattini.
\newblock Some remarks on special {K}\"ahler geometry.
\newblock {\em J. Geom. Phys.}, 59(7):755--763, 2009.

\bibitem[Bot95]{bottacin}
F.~Bottacin.
\newblock Symplectic geometry on moduli spaces of stable pairs.
\newblock {\em Ann. Sci. \'Ecole Norm. Sup. (4)}, 28(4):391--433, 1995.

\bibitem[BR94]{Biswas-Ramanan}
I.~Biswas and S.~Ramanan.
\newblock An infinitesimal study of the moduli of {H}itchin pairs.
\newblock {\em J. London Math. Soc. (2)}, 49(2):219--231, 1994.

\bibitem[BS02]{balaji_seshadri_1}
V.~Balaji and C.~S. Seshadri.
\newblock Semistable principal bundles. {I}. {C}haracteristic zero.
\newblock {\em J. Algebra}, 258(1):321--347, 2002.
\newblock Special issue in celebration of Claudio Procesi's 60th birthday.

\bibitem[CDDP15]{cddp}
W.~Chuang, D.-E. Diaconescu, R.~Donagi, and T.~Pantev.
\newblock Parabolic refined invariants and macdonald polynomials.
\newblock {\em Commun. Math. Phys.}, 335(3):1323--1379, 2015.

\bibitem[Che55]{chevalley}
C.~Chevalley.
\newblock Invariants of finite groups generated by reflections.
\newblock {\em Am. J. Math.}, 77(4):778--782, 1955.

\bibitem[DDD{\etalchar{+}}06]{geom_trans}
D.-E. Diaconescu, R.~Dijkgraaf, R.~Donagi, C.~Hofman, and T.~Pantev.
\newblock Geometric transitions and integrable systems.
\newblock {\em Nuclear Phys. B}, 752(3):329--390, 2006.

\bibitem[DDP07]{ddp}
D.-E. Diaconescu, R.~Donagi, and T.~Pantev.
\newblock Intermediate {J}acobians and {A}{D}{E} {H}itchin systems.
\newblock {\em Math. Res. Lett}, 14(5):745--756, 2007.

\bibitem[Del73]{Deligne_propre}
P.~Deligne.
\newblock Cohomologie a supports propres.
\newblock In {\em Th\'eorie des {T}opos et {C}ohomologie {E}tale des
  {S}ch\'emas}, Lecture Notes in Mathematics, Vol. 305, pages 250--480.
  Springer-Verlag, 1973.
\newblock S{\'e}minaire de G{\'e}om{\'e}trie Alg{\'e}brique du Bois-Marie
  1963--1964 (SGA 4).

\bibitem[Den95]{deninger_sign}
C.~Deninger.
\newblock Higher order operations in {D}eligne cohomology.
\newblock {\em Invent. Math.}, 120(2):289--315, 1995.

\bibitem[DG02]{don-gaits}
R.~Donagi and D.~Gaitsgory.
\newblock The gerbe of {H}iggs bundles.
\newblock {\em Transform. Groups}, 7(2):109--153, 2002.

\bibitem[DM93]{donagi_markman_cubic}
R.~Donagi and E.~Markman.
\newblock Cubics, integrable systems, and {C}alabi-{Y}au threefolds.
\newblock In {\em Proceedings of the Hirzebruch 65 conference on algebraic
  geometry (Ramat Gan)}, pages 199--221. 1993.

\bibitem[DM96]{donagi_markman}
R.~Donagi and E.~Markman.
\newblock Spectral curves, algebraically completely integrable, {H}amiltonian
  systems, and moduli of bundles.
\newblock In {\em Integrable systems and quantum groups ({M}ontecatini {T}erme,
  1993)}, volume 1620 of {\em L N M}, pages 1--119. Springer, Berlin, 1996.

\bibitem[Don93]{donagi_decomposition}
R.~Donagi.
\newblock Decomposition of spectral covers.
\newblock {\em Ast\'erisque}, 218:145--175, 1993.

\bibitem[Don95]{donagi_spectral_covers}
R.~Donagi.
\newblock Spectral {C}overs.
\newblock volume~28 of {\em M{S}{R}{I} {S}eries}. 1995.

\bibitem[Don98]{donagi_swis}
R.~Donagi.
\newblock Seiberg-{W}itten integrable systems.
\newblock In {\em Surveys in differential geometry: integrable systems}, Surv.
  Diff. Geom., IV, pages 83--129. Int. Press, Boston, MA, 1998.

\bibitem[DP05]{arijit_parthasarathi_hn}
Arijit Dey and R.~Parthasarathi.
\newblock On {H}arder-{N}arasimhan reductions for {H}iggs principal bundles.
\newblock {\em Proc. Indian Acad. Sci. Math. Sci.}, 115(2):127--146, 2005.

\bibitem[DP12]{don-pan}
R.~Donagi and T.~Pantev.
\newblock Langlands duality for {H}itchin systems.
\newblock {\em Invent. Math.}, 189(3):653--735, 2012.

\bibitem[DS95]{drinfeld_simpson}
V.G. Drinfeld and C.T. Simpson.
\newblock {$B$}-structures on {$G$}-bundles and local triviality.
\newblock {\em Math. Res. Lett.}, 2(6):823--829, 1995.

\bibitem[DW96]{donagi_witten}
R.~Donagi and E.~Witten.
\newblock Sypersymmetric {Y}ang-{M}ills theory and integrable systems.
\newblock {\em Nucl.Phys.B}, (460):299--334, 1996.

\bibitem[EM09]{ein_mustata_jets}
L.~Ein and M.~Musta{\c{t}}{\u{a}}.
\newblock Jet schemes and singularities.
\newblock In {\em Algebraic geometry---{S}eattle 2005. {P}art 2}, volume~80 of
  {\em Proc. Sympos. Pure Math.}, pages 505--546. Amer. Math. Soc., Providence,
  RI, 2009.

\bibitem[Fal93]{faltings}
G.~Faltings.
\newblock Stable {G}-bundles and projective connections.
\newblock {\em J. Algebraic Geom.}, 2:507--568, 1993.

\bibitem[FMW98]{friedman_morgan_witten}
R.~Friedman, J.~W. Morgan, and E.~Witten.
\newblock Principal {$G$}-bundles over elliptic curves.
\newblock {\em Math. Res. Lett.}, 5(1-2):97--118, 1998.

\bibitem[Fre99]{danfreed}
D.~S. Freed.
\newblock Special {K}\"ahler manifolds.
\newblock {\em Comm. Math. Phys.}, 203(1):31--52, 1999.

\bibitem[GO10]{bea_jh}
B.~Gra{\~n}a~Otero.
\newblock Jordan-{H}\"older reductions for principal {H}iggs bundles on curves.
\newblock {\em J. Geom. Phys.}, 60(11):1852--1859, 2010.

\bibitem[GPO14]{oscar_andre_compts}
O.~Garcia-Prada and A.~Oliveira.
\newblock Connectivity of {H}iggs bundle moduli for complex reductive {L}ie
  groups.
\newblock \href{http://lanl.arxiv.org/abs/1408.4778}{ar{X}iv 1408.4778}, 2014.

\bibitem[Gri68]{griffiths_periods}
Ph.~A. Griffiths.
\newblock Periods of integrals on algebraic manifolds, {I} and {I}{I}.
\newblock {\em American Journal of Mathematics}, 90(3):568--626, 805--865,
  1968.

\bibitem[Gro55]{Grothendieck_kansas}
A.~Grothendieck.
\newblock {\em A general theory of fibre spaces with structure sheaf}.
\newblock Univ. of Kansas, NSF Report, 1955.

\bibitem[Gro57]{Grothendieck_P1}
A.~Grothendieck.
\newblock Sur la classification des fibres holomorphes sur la sphere de
  {R}iemann.
\newblock {\em American Journal of Mathematics}, 79(1):121--138, 1957.

\bibitem[Her03]{hertling_ttstar}
C.~Hertling.
\newblock {$tt^*$} geometry, {F}robenius manifolds, their connections, and the
  construction for singularities.
\newblock {\em J. Reine Angew. Math.}, 555:77--161, 2003.

\bibitem[HHP10]{hertling_hoev_posthum}
C.~Hertling, L.~Hoevenaars, and H.~Posthuma.
\newblock Frobenius manifolds, projective special geometry and {H}itchin
  systems.
\newblock {\em J. Reine Angew. Math.}, 649:117--165, 2010.

\bibitem[Hit87a]{hitchin_sd}
N.~Hitchin.
\newblock The self-duality equations on a {R}iemann surface.
\newblock {\em Proc. London Math. Soc.}, 3(55):59--126, 1987.

\bibitem[Hit87b]{hitchin_sb}
N.~Hitchin.
\newblock Stable bundles and integrable systems.
\newblock {\em Duke Math. J.}, 54(1):91--114, 1987.

\bibitem[Hit92]{hitchin_teich}
N.~Hitchin.
\newblock Lie groups and {T}eichm\"uller space.
\newblock {\em Topology}, 31(3):449--473, 1992.

\bibitem[HM98]{hurtubise_markman_rk2}
J.~C. Hurtubise and E.~Markman.
\newblock Rank {$2$}-integrable systems of {P}rym varieties.
\newblock {\em Adv. Theor. Math. Phys.}, 2(3):633--695, 1998.

\bibitem[HM04]{hyeon_murphy}
D.~Hyeon and D.~Murphy.
\newblock Note on the stability of principal bundles.
\newblock {\em Proc. of AMS}, 132(8):2205--2213, 2004.

\bibitem[Hoe08]{hoevenaars}
L.~K. Hoevenaars.
\newblock Associativity for the {N}eumann system.
\newblock In {\em From {H}odge theory to integrability and {TQFT}
  tt*-geometry}, volume~78 of {\em Proc. Sympos. Pure Math.}, pages 215--238.
  AMS, Providence, RI, 2008.

\bibitem[Hur97]{hurtubise_kk}
J.~C. Hurtubise.
\newblock The algebraic geometry of the {K}ostant-{K}irillov form.
\newblock {\em J. London Math. Soc. (2)}, 56(3):504--518, 1997.

\bibitem[Huy06]{huybrechts_fm}
D.~Huybrechts.
\newblock {\em Fourier-{M}ukai transforms in algebraic geometry}.
\newblock Oxford Mathematical Monographs. The Clarendon Press Oxford University
  Press, Oxford, 2006.

\bibitem[Kji00]{kjiri}
M.~Kjiri.
\newblock The {$G$}-generalized {H}itchin systems and {P}rym varieties.
\newblock {\em J. Math. Phys.}, 41(11):7797--7807, 2000.

\bibitem[KKP08]{kkp_hodge_ms}
L~Katzarkov, M~Kontsevich, and T~Pantev.
\newblock Hodge theoretic aspects of mirror symmetry.
\newblock In {\em From Hodge theory to integrability and TQFT tt*-geometry.
  International workshop From TQFT to tt* and integrability, Augsburg, Germany,
  May 25--29, 2007}, pages 87--174. Providence, RI: American Mathematical
  Society (AMS), 2008.

\bibitem[Kod86]{kodaira_defo}
K.~Kodaira.
\newblock {\em Complex manifolds and deformation of complex structures}, volume
  283 of {\em Grundlehren der Mathematischen Wissenschaften [Fundamental
  Principles of Mathematical Sciences]}.
\newblock Springer-Verlag, New York, 1986.
\newblock Translated from the Japanese by Kazuo Akao, With an appendix by
  Daisuke Fujiwara.

\bibitem[KS58]{kodaira-spencer}
K.~Kodaira and D.~C. Spencer.
\newblock On deformations of complex analytic structures {I},{II}.
\newblock {\em Ann. of Math. (2)}, 67:328--466, 1958.

\bibitem[Las98]{laszlo_ell}
Y.~Laszlo.
\newblock About {G}-bundles over elliptic curves.
\newblock {\em Ann. Inst. Fourier}, 48:413--424, 1998.

\bibitem[Mar94]{markman_thesis}
E.~Markman.
\newblock Spectral curves and integrable systems.
\newblock {\em Compositio Math.}, 93(3):255--290, 1994.

\bibitem[Mar00]{markman_sw}
E.~Markman.
\newblock Algebraic geometry, integrable systems, and {S}eiberg-{W}itten
  theory.
\newblock In {\em Integrability: the {S}eiberg-{W}itten and {W}hitham equations
  ({E}dinburgh, 1998)}, pages 23--41. 2000.

\bibitem[Nit91]{nitsure}
N.~Nitsure.
\newblock Moduli space of semistable pairs on a curve.
\newblock {\em Proc. London Math. Soc. (3)}, 62(2):275--300, 1991.

\bibitem[Ram75]{ramanathan}
A.~Ramanathan.
\newblock Stable principal bundles on a compact {R}iemann surface.
\newblock {\em Math. Ann.}, 213:129--152, 1975.

\bibitem[Ram96a]{rama1}
A.~Ramanathan.
\newblock Moduli for principal bundles over algebraic curves. {I}.
\newblock {\em Proc. Indian Acad. Sci. Math. Sci.}, 106(3):301--328, 1996.

\bibitem[Ram96b]{rama}
A.~Ramanathan.
\newblock Moduli for principal bundles over algebraic curves. {II}.
\newblock {\em Proc. Indian Acad. Sci. Math. Sci.}, 106(4):421--449, 1996.

\bibitem[Sab11]{sabbah_nch}
C.~Sabbah.
\newblock Non-commutative {H}odge structures.
\newblock {\em Ann. Inst. Fourier}, 61(7):2681--2717, 2011.

\bibitem[Sco98]{scognamillo_elem}
R.~Scognamillo.
\newblock An elementary approach to the abelianization of the {H}itchin system
  for arbitrary reductive groups.
\newblock {\em Compositio Math.}, 110:17--37, 1998.

\bibitem[Ser58]{serre_fibres_algebriques}
J.-P. Serre.
\newblock Espaces fibr\'es alg\'ebriques.
\newblock In {\em S\'eminaire {C}laude {C}hevalley, {T}ome.\ 3, Exp.\ No.\ 1},
  pages 1--37. 1958.

\bibitem[Ses82]{seshadri_asterisque}
C.~S. Seshadri.
\newblock {\em Fibr\'es vectoriels sur les courbes alg\'ebriques}, volume~96 of
  {\em Ast\'erisque}.
\newblock Soci\'et\'e Math\'ematique de France, Paris, 1982.
\newblock Notes written by J.-M. Drezet from a course at the {\'E}cole Normale
  Sup{\'e}rieure, June 1980.

\bibitem[Sim88]{simpson_uniformisation}
C.~T. Simpson.
\newblock Constructing variations of {H}odge structure using {Y}ang-{M}ills
  theory and applications to uniformization.
\newblock {\em J. Amer. Math. Soc.}, 1(4):867--918, 1988.

\bibitem[Sim92]{hbls}
C.~Simpson.
\newblock Higgs bundles and local systems.
\newblock {\em Inst. {H}autes {\'E}tudes {S}ci. {P}ubl. {M}ath.}, (75):5--95,
  1992.

\bibitem[Sim94]{moduli2}
C.~Simpson.
\newblock Moduli of representations of the fundamental group of a smooth
  projective variety, {II}.
\newblock {\em Inst. {H}autes {\'E}tudes {S}ci. {P}ubl. {M}ath.}, (80):5--79,
  1994.

\bibitem[Ste65]{steinberg_regular}
R.~Steinberg.
\newblock Regular elements of semisimple algebraic groups.
\newblock {\em Inst. Hautes \'Etudes Sci. Publ. Math.}, (25):49--80, 1965.

\bibitem[Sze04]{szendroi_artin}
B.~Szendr{\H{o}}i.
\newblock Artin group actions on derived categories of threefolds.
\newblock {\em J. Reine Angew. Math.}, 572:139--166, 2004.

\bibitem[Sze08]{szendroi_3folds}
B.~Szendr{\H{o}}i.
\newblock Sheaves on fibered threefolds and quiver sheaves.
\newblock {\em Comm. Math. Phys.}, 278(3):627--641, 2008.

\bibitem[Tu93]{tu_elliptic}
L.~Tu.
\newblock Semistable bundles over an elliptic curve.
\newblock {\em Advances in Math.}, 98:1--26, 1993.

\bibitem[Wei83]{weinstein_local_structure}
A.~Weinstein.
\newblock The local structure of {P}oisson manifolds.
\newblock {\em J. Differential Geom.}, 18(3):523--557, 1983.

\end{thebibliography}
\end{document}